\documentclass[11pt]{article}
\usepackage{amsfonts}
\usepackage{color}
\newtheorem{theo}{Theorem}[section]
\newtheorem{lem}[theo]{Lemma}
\newtheorem{cor}[theo]{Corollary}

\newtheorem{defi}{Definition}[section]

\newcommand{\mysection}[1]{\section{#1} \setcounter{equation}{0}}
\newcommand{\proof}{{\sc Proof.} \quad}
\newcommand{\proofc}{{\sc Proof} \ }
\newcommand{\be}{\begin{equation} \label}
\newcommand{\ee}{\end{equation}}
\newcommand{\bea}{\begin{eqnarray}\label}
\newcommand{\eea}{\end{eqnarray}}
\newcommand{\bas}{\begin{eqnarray*}}
\newcommand{\eas}{\end{eqnarray*}}
\newcommand{\bit}{\begin{itemize}}
\newcommand{\eit}{\end{itemize}}
\newcommand{\qed}{\hfill$\Box$ \vskip.2cm}
\newcommand{\nn}{\nonumber}
\newcommand{\R}{\mathbb{R}}
\newcommand{\N}{\mathbb{N}}
\newcommand{\pO}{\partial\Omega}

\newcommand{\eps}{\varepsilon}

\newcommand{\wto}{\rightharpoonup}
\newcommand{\wsto}{\stackrel{\star}{\rightharpoonup}}
\newcommand{\hra}{\hookrightarrow}
\newcommand{\io}{\int_\Omega}
\newcommand{\bom}{\overline{\Omega}}

\newcommand{\abs}{\\[5pt]}
\newcommand{\Abs}{\\[5mm]}

\newcommand{\proj}{{\cal P}}
\newcommand{\neps}{n_\eps}
\newcommand{\ceps}{c_\eps}
\newcommand{\ueps}{u_\eps}
\newcommand{\Peps}{P_\eps}
\newcommand{\tme}{T_{max,\eps}}
\newcommand{\htme}{\widehat{T}_{max,\eps}}
\newcommand{\kd}{k_D}
\newcommand{\onz}{\overline{n_0}}
\newcommand{\ts}{t_\star}
\newcommand{\ps}{p_\star}
\newcommand{\pss}{p^\star}
\newcommand{\tpsi}{\widetilde{\psi}}
\newcommand{\hc}{\widehat{c}}
\newcommand{\hu}{\widehat{u}}
\begin{document}
\enlargethispage{10mm}
\title{Global existence and stabilization in a degenerate chemotaxis-Stokes system with 
mildly strong diffusion enhancement}
\author{
Michael Winkler\footnote{michael.winkler@math.uni-paderborn.de}\\
{\small Institut f\"ur Mathematik, Universit\"at Paderborn,}\\
{\small 33098 Paderborn, Germany} }
\date{}
\maketitle
\begin{abstract}
\noindent 
  A class of chemotaxis-Stokes systems generalizing the prototype
  \bas
    	\left\{ \begin{array}{rcl}
    	n_t + u\cdot\nabla n &=& \nabla \cdot \big(n^{m-1}\nabla n\big)  - \nabla \cdot \big(n\nabla c\big), \\[1mm]
    	c_t + u\cdot\nabla c &=& \Delta c-nc, \\[1mm]
     	u_t +\nabla P  &=& \Delta u + n \nabla \phi, \qquad \nabla\cdot u =0, 
    	\end{array} \right.
  \eas
  is considered in bounded convex three-dimensional domains, where $\phi\in W^{2,\infty}(\Omega)$ is given.\abs
  The paper develops an analytical approach which consists in a combination of energy-based arguments and
  maximal Sobolev regularity theory, and which allows for the construction of global bounded weak solutions
  to an associated initial-boundary value problem under the assumption that 
  \be{m0}
	m>\frac{9}{8}.
  \ee
  Moreover, the obtained solutions are shown to approach the spatially homogeneous steady state
  $(\frac{1}{|\Omega|} \io n_0,0,0)$ in the large time limit.\abs
  This extends previous results which either relied on different and apparently less significant 
  energy-type structures, or on completely alternative approaches, and thereby exclusively achieved 
  comparable results under hypotheses stronger than (\ref{m0}).\abs
\end{abstract}
\newpage
\section{Introduction}\label{intro}
We consider the chemotaxis-Stokes system
\be{0}
    	\left\{ \begin{array}{rcll}
    	n_t + u\cdot\nabla n &=& \nabla \cdot \Big(D(n)\nabla n\Big)  - \nabla \cdot \big(n\nabla c\big), 
	\qquad & x\in\Omega, \ t>0,\\[1mm]
    	c_t + u\cdot\nabla c &=& \Delta c-nc, 
	\qquad & x\in\Omega, \ t>0,   \\[1mm]
     	u_t +\nabla P  &=& \Delta u + n \nabla \phi, \qquad \nabla\cdot u =0, 
	\qquad & x\in\Omega, \ t>0, 
    	\end{array} \right.
\ee
which was proposed in \cite{goldstein2005} and \cite{DiFLM} as a model for the spatio-temporal evolution
in populations of oxytactically moving bacteria that interact with a surrounding fluid
through transport and buoyancy, where $n, c, u$ and $P$ denote the density of cells, the oxygen concentration,
the fluid velocity and its associated pressure, respectively, and where the diffusivity $D$ and the
gravitational potential $\phi=\phi(x)$ are given smooth parameter functions
(cf.~also \cite{bellomo} for a recent independent derivation of (\ref{0})
on the basis of fundamental principles from the kinetic theory of active particles).
Indeed, as reported in \cite{goldstein2004} and \cite{goldstein2005},
even in such a simple setting lacking any reinforcement of chemotactic motion by signal production 
through cells, quite a colorful collective behavior can be observed, including the formation of aggregates
and the emergence of large-scale convection patterns.\abs
In modification of the original model from \cite{goldstein2005} in which $D\equiv 1$, the authors in \cite{DiFLM}
suggested to adequately account for the finite size of bacteria by assuming that the random movement of cells
is nonlinearly enhanced at large densities, leading to the choice 
\be{proto}
	D(s)=s^{m-1}
	\qquad \mbox{for } s\ge 0
\ee
with some $m>1$
in the prototypical case of porous medium type diffusion. 
In comparison to the case $D\equiv 1$, nonlinear diffusion mechanisms of this type may suppress the occurrence of
blow-up phenomena, as known to be enforced by chemotactic cross-diffusion 
e.g.~in frameworks such as that 
addressed by the classical Keller-Segel system (\cite{herrero_velazquez}, \cite{win_JMPA}).
In fact, 
in three-dimensional initial value problems for (\ref{0}) with $D\equiv 1$, global smooth and bounded
solutions could be shown to exist only under appropriate smallness assumptions on the initial data
(\cite{DLM}, \cite{kozono_miura_sugiyama}, \cite{chae_kang_lee_CPDE}, \cite{cao_lankeit}),
while for arbitrarily large data so far only certain global weak 
solutions have been constructed, which do become smooth eventually but may develop singularities
prior to such ultimate regularization (\cite{win_CPDE}, \cite{win_TRAN}).
Contrary to this, assuming (\ref{proto}) to hold,
recent analysis has revealed the condition
\be{m1}
	m>\frac{7}{6}
\ee
as sufficient for global existence and boundedness of weak solutions to an associated no-flux-no-flux-Dirichlet
initial-boundary value problem for all reasonably regular initial data
in three-dimensional bounded convex domains (\cite{win_CVPDE}, cf.~also \cite{liu_lorz}).
This partially extended a precedent result which asserted global solvability within the larger range
$m>\frac{8}{7}$,
but only in a class of weak solutions locally bounded in $\bom\times [0,\infty)$ (\cite{taowin_ANIHPC}).
For smaller values of $m>1$, up to now existence results are limited to classes of possibly unbounded
solutions (\cite{duan_xiang_IMRN}).\abs
In view of lacking complementary results on possibly occurring singularity formation phenomena, the question of 
identifying an {\em optimal} condition on $m\ge 1$ ensuring global boundedness in the three-dimensional 
version of (\ref{0}) remains an open challenge,
thus marking a substantial difference to the two-dimensional situation 
in which global existence and boundedness results 
are available for several variants of (\ref{0}) already in presence of linear cell diffusion, and even when
the fluid flow is governed by the corresponding full nonlinear Navier-Stokes system
(\cite{DLM}, \cite{win_CPDE}, \cite{win_ARMA}, \cite{chae_kang_lee2015}, \cite{zhang_li}).\Abs
{\bf Main results.} \quad
It is the purpose of this work to demonstrate how an adequate combination of energy-based arguments and 
maximal Sobolev regularity theory can be used to further advance the analysis of (\ref{0}), with $D$ essentially
of the form in (\ref{proto}), even in previously unexplored ranges of $m$.
In fact, in the first step our approach we will make use of an observation to be stated in Lemma \ref{lem1},
according to which the system (\ref{0}) also for $m>1$ 
continues to feature an energy-type structure known to be present when $m=1$ even in an
associated chemotaxis-Navier-Stokes system (\cite{win_ANIHPC}; cf.~also \cite{DLM} and \cite{win_CPDE} for precedent
partial findings in this direction).
By means of a first iterative bootstrap procedure, 
the correspondingly obtained a priori estimates will be turned into some regularity information on the 
solution component $n$ (Section \ref{sect4} and Section \ref{sect5}), 
which itself can be used as a starting point for a second recursive argument:
Namely, investigating how far regularity information of the latter type influences integrability properties
of $u$ and $\nabla c$ through maximal Sobolev regularity estimates (Section \ref{sect6}), we will be able to successively 
improve our knowledge on available integral bounds for all solution components
under the mild assumption that in the setup of (\ref{proto}) we merely have
\be{m2}
	m>\frac{9}{8}
\ee
(Section \ref{sect7} and Section \ref{sect8}).
The estimates thereby obtained will provide appropriate compactness properties which will firstly allow us
to construct global bounded weak solutions to (\ref{0}) via a suitable approximation procedure (Section \ref{sect9}), 
and which thereafter secondly enable us to assert stabilization toward spatially homogeneous equilibria
(Section \ref{sect10}).\abs
In order to formulate our results in these directions, let us specify the setup of our analysis by declaring that
throughout the sequel we shall assume $D$ to generalize the choice in (\ref{proto}) in that
\be{D}
	D\in C^{\vartheta}_{loc}([0,\infty)) \cap C^2((0,\infty))
	\quad \mbox{is such that} \quad
	D(s) \ge \kd s^{m-1}
	\quad \mbox{ for all $s\ge 0$}
\ee
with some $\vartheta\in (0,1), \kd>0$ and $m>1$,
and by considering the initial-boundary value problem for (\ref{0}) associated with the requirements that
\be{0i}
    n(x,0)=n_0(x), \quad c(x,0)=c_0(x) \quad \mbox{and} \quad u(x,0)=u_0(x), \qquad x\in\Omega,
\ee
as well as
\be{0b}
    	\Big(D(n)\nabla n - n\nabla c\Big) \cdot \nu=0, \quad \frac{\partial c}{\partial\nu}=0 
	\quad \mbox{and} \quad u=0
	\qquad \mbox{on } \pO,
\ee
in a bounded convex domain $\Omega\subset \R^3$ with smooth boundary.
As for the initial data herein, we shall suppose for convenience that
\be{init}
    \left\{
    \begin{array}{l}
    n_0 \in C^\omega(\bom) \quad \mbox{for some $\omega>0$ with } n_0\ge 0 \mbox{ in $\Omega$ and }
	n_0 \not\equiv 0, \quad
	\mbox{that} \\
    c_0 \in W^{1,\infty}(\Omega) \quad \mbox{ satisfies $c_0\ge 0$ in $\Omega$, \quad and that}\\
    u_0 \in D(A^{\alpha}) \quad \mbox{for some $\alpha\in (\frac{3}{4},1)$,}
    \end{array}
    \right.
\ee
where $A$ denotes the Stokes operator in $L^2_\sigma(\Omega):=\{ \varphi\in L^2(\Omega) \ | \ \nabla \cdot \varphi=0 \}$
with its domain given by $D(A):=W^{2,2}(\Omega)\cap W_0^{1,2}(\Omega) \cap L^2_\sigma(\Omega)$ (\cite{sohr_book}).\abs
We shall then obtain the following result on global existence and large time behavior,
where as in several places below we make use of the abbreviation 
$\overline{\varphi}:=\frac{1}{|\Omega|} \io \varphi$ for $\varphi\in L^1(\Omega)$.
\begin{theo}\label{theo99}
  Let $\Omega\subset\R^3$ be a bounded convex domain with smooth boundary and 
  $\phi\in W^{2,\infty}(\Omega)$, and suppose that $D$ is such that 
  (\ref{D}) holds with some
  \be{m}
	m>\frac{9}{8}.
  \ee
  Then for each $n_0, c_0$ and $u_0$ satisfying (\ref{init}) there exist functions
  \be{99.1}
	\left\{ \begin{array}{l}
	n \in L^\infty(\Omega\times (0,\infty)) \cap C^0([0,\infty);(W_0^{2,2}(\Omega))^\star), \\
	c \in \bigcap_{p>1} L^\infty((0,\infty); W^{1,p}(\Omega)) 
	\cap C^0(\bom\times [0,\infty)) \cap C^{1,0}(\bom\times (0,\infty)), \\
	u\in L^\infty(\Omega\times (0,\infty)) \cap L^2_{loc}([0,\infty);W_0^{1,2}(\Omega)\cap L^2_\sigma(\Omega))
	\cap C^0(\bom\times [0,\infty))
	\end{array} \right.
  \ee
  such that the triple $(n,c,u)$ forms a global weak solution of (\ref{0}), (\ref{0i}), (\ref{0b})
  in the sense of Definition \ref{defi_weak} below.\\
  Moreover, this solution has the property that for arbitrary $p\ge 1$ we have
  \be{99.2}
    	\|n(\cdot,t)-\onz\|_{L^p(\Omega)} 
	+ \|c(\cdot,t)\|_{W^{1,\infty}(\Omega)}
	+ \|u(\cdot,t)\|_{L^\infty(\Omega)}
	\to 0
	\qquad \mbox{as } t\to\infty.
  \ee
\end{theo}
As a by-product, this trivially extends previous results on blow-up suppression in the associated fluid-free
chemotaxis system with porous medium-type diffusion and signal consumption, as obtained on letting 
$u\equiv 0$ in (\ref{0}). Even for the latter, apparently somewhat simpler system, 
only under the assumption (\ref{m1}) global bounded solutions have been known to 
exist (\cite{liangchen_wang_et_al_ZAMP2016}), with again no example of blow-up available for any choice of
$D$ yet.\abs
In order to further put these results in perspective, let us note that alternative modeling approaches
suggest to introduce as blow-up inhibiting mechanisms certain saturation effects in the cross-diffusive term 
in (\ref{0}) at large cell densities (cf.~e.g.~the survey \cite{hillen_painter2009}).
Indeed, if in (\ref{0}) the summand $-\nabla\cdot (n\nabla c)$ is replaced by 
$-\nabla\cdot (nS(n)\nabla c)$ with $S$ suitably generalizing
the prototype given by $S(s)=(s+1)^{-\alpha}$ for all $s\ge 0$ and some $\alpha>0$, then known results assert
global existence of bounded solutions to a corresponding initial-boundary value problem
when in the context of (\ref{D}) we have $m+\alpha>\frac{7}{6}$ (\cite{yilong_wang_xie}), 
which in the particular case $\alpha=0$
considered here rediscovers (\ref{m1}) and is thereby stronger than (\ref{m}).
An interesting open problem, partially addressed in \cite{cao_wang}, \cite{wangyulan_CAMWA} and \cite{wang_xiang},
consists in determining optimal conditions on the interplay between these two mechanisms which indeed prevent explosions.
\mysection{Approximation by non-degenerate problems}\label{sect_approx}
In order to construct solutions of (\ref{0}) through an appropriate approximation,
following natural regularization procedures we fix a family $(D_\eps)_{\eps\in (0,1)}$ of functions
\bea{Deps}
	& & D_\eps \in C^2([0,\infty))
	\quad \mbox{such that} \quad
	D_\eps(s) \ge \eps \quad \mbox{for all $s\ge 0$ and } \eps\in (0,1)
	\quad \mbox{and} \quad \nn\\
	& & D(s) \le D_\eps(s) \le D(s)+2\eps
	\quad \mbox{for all $s\ge 0$ and } \eps\in (0,1),
\eea
and we moreover regularize the cross-diffusive term in (\ref{0}) by introducing a family
$(\chi_\eps)_{\eps\in (0,1)} \subset C_0^\infty([0,\infty))$ fulfilling
\be{chi}
	0 \le \chi_\eps \le 1 \mbox{ in $[0,\infty)$,}
	\quad
	\chi_\eps\equiv 1 \mbox{ in $[0,\frac{1}{\eps}]$ \quad and} \quad
	\chi_\eps\equiv 0 \mbox{ in $[\frac{2}{\eps},\infty)$,}
\ee
and by letting
\be{def_F}
	F_\eps(s):=\int_0^s \chi_\eps(\sigma) d\sigma, \qquad s\ge 0,
\ee
for $\eps\in (0,1)$. Then $F_\eps \in C^\infty([0,\infty))$ satisfies
\be{F}
	0 \le F_\eps(s) \le s
	\quad \mbox{and} \quad
	0 \le F_\eps'(s) \le 1	
	\qquad \mbox{for all } s\ge 0
\ee
as well as
\be{F_conv}
	F_\eps(s)\nearrow s
	\quad \mbox{for all } s\ge 0
	\quad \mbox{and} \quad
	F_\eps'(s) \nearrow 1
	\quad \mbox{for all } s>0
	\qquad \mbox{as $\eps\searrow 0$}.
\ee
These choices in particular guarantee that each of the approximate variants of (\ref{0}), (\ref{0i}), (\ref{0b})
given by
\be{0eps}
	\left\{ \begin{array}{ll}
	\partial_t n_\eps + \ueps \cdot \nabla \neps = \nabla \cdot \Big(D_\eps(\neps)\nabla \neps \Big) 
	- \nabla \cdot \Big(\neps F_\eps'(\neps) \nabla \ceps \Big),
	\qquad & x\in \Omega, \ t>0, \\[1mm]
	\partial_t c_\eps+ \ueps \cdot \nabla \ceps = \Delta \ceps - F_\eps(\neps)\ceps,
	\qquad & x\in \Omega, \ t>0, \\[1mm]
	\partial_t u_\eps + \nabla \Peps = \Delta \ueps + \neps \nabla \phi,
	\qquad & x\in \Omega, \ t>0, \\[1mm]
	\nabla \cdot \ueps=0,
	\qquad & x\in \Omega, \ t>0, \\[1mm]
	\frac{\partial \neps}{\partial\nu} = \frac{\partial \ceps}{\partial\nu}=0, \quad \ueps=0,
	\qquad & x\in \partial\Omega, \ t>0, \\[1mm]
	\neps(x,0)=n_0(x), \quad \ceps(x,0)=c_0(x), \quad \ueps(x,0)=u_0(x),
	\qquad & x\in\Omega,
 	\end{array} \right.
\ee
for $\eps\in (0,1)$, possesses globally defined classical solutions:
\begin{lem}\label{lem0}
  Let $\eps\in (0,1)$. Then there exist functions
  \bas
	\left\{ \begin{array}{l}
	\neps\in C^0(\bom \times [0,\infty)) \cap C^{2,1}(\bom\times (0,\infty)), \\
	\ceps\in C^0(\bom \times [0,\infty)) \cap C^{2,1}(\bom\times (0,\infty)), \\
	\ueps\in C^0(\bom \times [0,\infty)) \cap C^{2,1}(\bom\times (0,\infty)), \\
	\Peps\in C^{1,0}(\bom\times (0,\infty)),
	\end{array} \right.
  \eas
  such that $(\neps,\ceps,\ueps,\Peps)$ solves (\ref{0eps}) classically in $\Omega \times (0,\infty)$, and such
  that $\neps$ and $\ceps$ are nonnegative in $\Omega\times (0,\infty)$.
\end{lem}
\proof
  By means of standard arguments from the local existence theories of taxis-type cross diffusive parabolic systems
  and the Stokes evolution equation (\cite{amann}, \cite{sohr_book}, \cite{lankeit_consumption}, \cite{win_CPDE}), 
  it follows that there exist
  $\tme\in (0,\infty]$ and at least one classical solution $(\neps,\ceps,\ueps,P_\eps) \in 
  \Big(C^0(\bom\times [0,\tme);\R^5) \cap C^{2,1}(\bom\times (0,\tme);\R^5)\Big) \times C^{1,0}(\bom\times (0,\tme))$
  which is such that $\neps\ge 0$ and $\ceps\ge 0$ in $\Omega\times (0,\tme)$,
  that $\ceps\in C^0([0,\tme);W^{1,p}(\Omega))$ for all $p\ge 1$ and that if $\tme<\infty$ then
  \be{0.1}
	\limsup_{t\nearrow\tme} \Big( \|\neps(\cdot,t)\|_{C^2(\bom)}
	+ \|\ceps(\cdot,t)\|_{C^2(\bom)} + \|\ueps(\cdot,t)\|_{C^2(\bom)} \Big) = \infty.
  \ee
  For each $T>0$, however, using that for any fixed $\eps\in (0,1)$ the function $F_\eps'$ has its support
  located in $[0,\frac{2}{\eps}]$ according to (\ref{def_F}) and (\ref{chi}), successive application of
  well-established $L^p$ estimation techniques and methods from higher order regularity theories for scalar
  parabolic equations and the Stokes system yields $C_1(\eps,T)>0$ such that
  \bas
	\|\neps(\cdot,t)\|_{C^2(\bom)}
	+ \|\ceps(\cdot,t)\|_{C^2(\bom)}
	+ \|\ueps(\cdot,t)\|_{C^2(\bom)}
	\le C_1(\eps,T)
	\qquad \mbox{for all } t\in (\tau_\eps,\htme),
  \eas
  where $\tau_\eps:=\min\{\frac{1}{2} T, \frac{1}{2} \tme\}$ and $\htme:=\min\{T,\tme\}$.
  This shows that (\ref{0.1}) cannot hold when $\tme$ is finite, whence we actually must have $\tme=\infty$.
\qed
In order to simplify presentation, throughout the sequel we shall tacitly assume that $(n_0,c_0,u_0)$ satisfies
(\ref{init}), and that for $\eps\in (0,1)$, $(\neps,\ceps,\ueps,\Peps)$ denotes the corresponding solution to 
(\ref{0eps}) obtained in Lemma \ref{lem0}.\abs
The following two basic properties thereof		
are immediate consequences of an integration in the first equation therein, as well as an application of the 
maximum principle to the second.
\begin{lem}\label{lem_basic}
  We have
  \be{mass}
	\|\neps(\cdot,t)\|_{L^1(\Omega)}=\io n_0 
	\qquad \mbox{for all } t>0
  \ee
  as well as
  \be{cinfty}
	\|\ceps(\cdot,t)\|_{L^\infty(\Omega)} \le \|c_0\|_{L^\infty(\Omega)}
	\qquad \mbox{for all } t>0.
  \ee
\end{lem}
\mysection{Directly exploiting the natural quasi-energy structure of (\ref{0eps})}\label{sect3}
Some first regularity properties beyond those from Lemma \ref{lem_basic} can be obtained by making use
of a quasi-energy structure which the approximate problems (\ref{0eps}) inherit from (\ref{0}) thanks to
the particular link between the dependence on $\neps$ of the interaction terms
$-\nabla\cdot (\neps F_\eps'(\neps)\nabla\ceps)$ and $-F_\eps(\neps)\ceps$ therein.
Similar energy-like properties have been used in previous studies on related problems (\cite{DLM}, \cite{taowin_ANIHPC},
\cite{win_CPDE}), but only in few cases the fluid velocity has been included (\cite{lankeit_M3AS}, \cite{win_ANIHPC},
\cite{win_TRAN}).
\begin{lem}\label{lem1}
  There exist $\kappa>0$ and $C>0$ such that
  \bea{1.1}
	& & \hspace*{-20mm}
	\frac{d}{dt} \bigg\{ \io \neps \ln \neps
	+ \frac{1}{2} \io \frac{|\nabla\ceps|^2}{\ceps}
	+ \kappa \io |\ueps|^2 \bigg\} \nn\\
	& & + \frac{1}{C} \bigg\{ \io \neps \ln \neps
	+ \frac{1}{2} \io \frac{|\nabla\ceps|^2}{\ceps}
	+ \kappa \io |\ueps|^2 \bigg\} \nn\\
	& & + \frac{1}{C} \bigg\{ \io \neps^{m-2} |\nabla\neps|^2
	+ \io \frac{|\nabla\ceps|^4}{\ceps^3}
	+ \io |\nabla\ueps|^2 \bigg\}
	\le C
	\qquad \mbox{for all } t>0.
  \eea
\end{lem}
\proof
  The derivation of (\ref{1.1}) follows a standard reasoning combining ideas from \cite{DLM}, \cite{win_CPDE}
  and \cite{win_ANIHPC}:
  By means of straightforward computation using the first two equations in (\ref{0eps}) (cf.~\cite[Lemma 3.2]{win_CPDE}
  for details), we obtain the identity
  \bea{1.2}
	& & \hspace*{-20mm}
	\frac{d}{dt} \bigg\{ \io \neps\ln \neps + \frac{1}{2} \io \frac{|\nabla\ceps|^2}{\ceps} \bigg\} 
	+ \io \frac{D_\eps(\neps)}{\neps} |\nabla\neps|^2
	+ \io \ceps |D^2 \ln \ceps|^2 \nn\\
	&=& - \frac{1}{2} \io \frac{|\nabla\ceps|^2}{\ceps^2} (\ueps\cdot\nabla\ceps)
	+ \io \frac{\Delta\ceps}{\ceps} (\ueps\cdot\nabla\ceps) \nn\\
	& & - \frac{1}{2} \io F_\eps(\neps) \frac{|\nabla\ceps|^2}{\ceps}
	+ \frac{1}{2} \int_{\pO} \frac{1}{\ceps} \frac{\partial |\nabla\ceps|^2}{\partial\nu}
	\qquad \mbox{for all } t>0,
  \eea
  where
  \be{1.22}
	\io \frac{D_\eps(\neps)}{\neps} |\nabla\neps|^2
	\ge \kd \io \neps^{m-2} |\nabla\neps|^2
	\qquad \mbox{for all } t>0
  \ee
  by (\ref{Deps}) and (\ref{D}), and where the two last summands on the right are nonpositive by nonnegativity
  of $F_\eps$ and due to the fact that $\frac{\partial |\nabla\ceps|^2}{\ceps}\le 0$ on
  $\pO\times (0,\infty)$ thanks to the convexity of $\Omega$ (\cite[Lemme 2.I.1]{lions_ARMA}).
  We next recall from \cite[Lemma 3.3]{win_CPDE} that
  \bas
	\io \frac{|\nabla\ceps|^4}{\ceps^3} \le C_1 \io \ceps |D^2 \ln \ceps|^2
	\qquad \mbox{for all } t>0
  \eas
  with $C_1:=(2+\sqrt{3})^2$, 
  and, after two integrations by parts in (\ref{1.2}), combine (\ref{cinfty}) with Young's inequality to estimate
  \bea{1.4}
	- \frac{1}{2} \io \frac{|\nabla\ceps|^2}{\ceps^2} (\ueps\cdot\nabla\ceps)
	+ \io \frac{\Delta\ceps}{\ceps} (\ueps\cdot\nabla\ceps)
	&=& \frac{1}{2} \io \frac{|\nabla\ceps|^2}{\ceps^2} (\ueps\cdot\nabla\ceps)
	- \io \frac{1}{\ceps} \nabla\ceps \cdot (D^2\ceps\cdot\nabla\ueps) \nn\\
	& & - \io \frac{1}{\ceps} \nabla\ceps \cdot (\nabla\ueps\cdot\nabla\ceps) \nn\\
	&=& 
	- \io \frac{1}{\ceps} \nabla\ceps \cdot (\nabla\ueps\cdot\nabla\ceps) \nn\\
	&\le& \frac{1}{2C_1} \io \frac{|\nabla \ceps|^4}{\ceps^3}
	+ C_2 \io |\nabla\ueps|^2
	\qquad \mbox{for all } t>0
  \eea
  with $C_2:=\frac{1}{2} \|c_0\|_{L^\infty(\Omega)}$.
  Now testing the third equation in (\ref{0eps}) by $\ueps$, thanks to the continuity of the embeddings
  $W_0^{1,2}(\Omega) \hra L^6(\Omega)$ and $W^{1,2}(\Omega) \hra L^\frac{12}{5m}(\Omega)$
  we independently see using the Gagliardo-Nirenberg inequality, Young's inequality and (\ref{mass}) that
  there exist positive constants $C_3, C_4, C_5$ and $C_6$ such that
  \bas
	\frac{1}{2} \frac{d}{dt} \io |\ueps|^2
	+ \io |\nabla\ueps|^2
	&=& \io \neps\ueps\cdot\nabla\phi \\
	&\le& \|\nabla\phi\|_{L^\infty(\Omega)} \|\ueps\|_{L^6(\Omega)} \|\neps\|_{L^\frac{6}{5}(\Omega)} \\
	&\le& C_3\|\nabla \ueps\|_{L^2(\Omega)} \|\neps^\frac{m}{2}\|_{L^\frac{12}{5m}(\Omega)}^\frac{2}{m} \\
	&\le& \frac{1}{2} \io |\nabla\ueps|^2
	+ \frac{C_3^2}{2} \|\neps^\frac{m}{2}\|_{L^\frac{12}{5m}(\Omega)}^\frac{4}{m} \\
	&\le& \frac{1}{2} \io |\nabla\ueps|^2
	+ C_4 \cdot \Big\{ \|\nabla\neps^\frac{m}{2}\|_{L^2(\Omega)}^\frac{2}{3m-1} 
		\|\neps^\frac{m}{2}\|_{L^\frac{2}{m}(\Omega)}^\frac{10m-4}{3m^2-m}
	+ \|\neps^\frac{m}{2}\|_{L^\frac{2}{m}(\Omega)}^\frac{4}{m} \Big\} \nn\\
	&\le& \frac{1}{2} \io |\nabla\ueps|^2
	+ C_5 \|\nabla\neps^\frac{m}{2}\|_{L^2(\Omega)}^\frac{2}{3m-1}+C_5 \\
	&\le& \frac{1}{2} \io |\nabla\ueps|^2
	+ \frac{\kd}{4(C_2+1)} \io \neps^{m-2} |\nabla\neps|^2 + C_6
	\qquad \mbox{for all } t>0.
  \eas
  In combination with (\ref{1.2}), (\ref{1.22}) and (\ref{1.4}), this shows that
  \bas
	& & \hspace*{-20mm}
	\frac{d}{dt} \bigg\{ \io \neps\ln\neps + \frac{1}{2} \io \frac{|\nabla\ceps|^2}{\ceps}
	+ (C_2+1) \io |\ueps|^2 \bigg\} \\
	& & + \frac{\kd}{2} \io \neps^{m-2} |\nabla\neps|^2
	+ \frac{1}{2C_1} \io \frac{|\nabla\ceps|^4}{\ceps^3}
	+ \io |\nabla\ueps|^2 \\
	&\le& 2(C_2+1)C_6
	\qquad \mbox{for all } t>0.
  \eas
  Since finally from the Gagliardo-Nirenberg inequality along with Young's inequality and (\ref{cinfty})
  we readily obtain $C_7>0$ such that
  \bas
	\io \neps\ln\neps + \frac{1}{2} \io \frac{|\nabla\ceps|^2}{\ceps} + (C_2+1) \io |\ueps|^2
	\le C_7 \cdot \bigg\{ \io \neps^{m-2} |\nabla\neps|^2
	+ \io \frac{|\nabla\ceps|^4}{\ceps^3} + \io |\nabla\ueps|^2 + 1 \bigg\}
  \eas
  for all $t>0$, this readily establishes (\ref{1.1}) upon evident obvious choices of $\kappa$ and $C$.
\qed
In the sequel we shall make use of the latter exclusively through the following direct consequences.
\begin{lem}\label{lem2}
  There exists $C>0$ such that for all $\eps\in (0,1)$,
  \be{2.1}
	\int_t^{t+1} \io |\nabla\neps^\frac{m}{2}|^2 \le C
	\qquad \mbox{for all } t\ge 0
  \ee
  and
  \be{2.2}
	\int_t^{t+1} \io |\nabla\ceps|^4 \le C
	\qquad \mbox{for all } t\ge 0
  \ee
  as well as
  \be{2.3}
	\int_t^{t+1} \io |\nabla\ueps|^2 \le C
	\qquad \mbox{for all } t\ge 0.
  \ee
\end{lem}
\proof
  All inequalities immediately result from an integration of (\ref{1.1}) because of (\ref{cinfty}) and the fact that
  $\io \neps\ln \neps \ge -\frac{|\Omega|}{e}$ for all $t\ge 0$.
\qed
\mysection{Preparing an inductive argument}\label{sect4}
We next address the question how far
an informational background such as the one provided by Lemma \ref{lem2} and Lemma \ref{lem_basic}
can be exploited so as to derive further regularity features of solutions to (\ref{0eps}).
More precisely, we shall be concerned with the problem of finding appropriate conditions on $m$ and the numbers
$\ps\ge 1$ and $\pss>\ps$ such that bounds of the form
\be{bound}
	\io \neps^{p}(\cdot,t) \le C
	\quad \mbox{and} \quad
	\int_t^{t+1} \io \Big|\nabla\neps^\frac{p+m-1}{2}\Big|^2
	\le C
	\qquad \mbox{for all } t\ge 0,
\ee
assumed to be present for $p=\ps$, 
can be shown to imply the same estimates for the corresponding quantities for $p=\pss$.\abs
Our first result in this direction actually requires a bound for $\neps$ in the single space
$L^\infty((0,\infty);L^{\ps}(\Omega))$ only, but additionally relies on a space-time regularity property
of $\nabla\ceps$ in asserting the following.
\begin{lem}\label{lem8}
  Let $m>1, \ps\ge 1, p>1$ and $q\ge 2$ be such that
  \be{8.1}
	p\le \frac{2(q-1)}{3} \ps + (2q-1)(m-1).
  \ee
  Then for all $K>0$ there exists $C=C(\ps,p,q,K)>0$ such that if for some $\eps\in (0,1)$ we have
  \be{8.2}
	\io \neps^{\ps}(\cdot,t) \le K
	\qquad \mbox{for all } t\ge 0
  \ee
  and
  \be{8.3}
	\int_t^{t+1} \io |\nabla\ceps|^{2q} \le K
	\qquad \mbox{for all } t\ge 0,
  \ee
  then
  \be{8.4}
	\io \neps^p(\cdot,t) \le C
	\qquad \mbox{for all } t\ge 0
  \ee
  and
  \be{8.5}
	\int_t^{t+1} \io \Big| \nabla \neps^\frac{p+m-1}{2}\Big|^2 \le C
	\qquad \mbox{for all } t\ge 0.
  \ee
\end{lem}
\proof
  In view of (\ref{mass}) and Lemma \ref{lem2}, 
  since $\frac{2(q-1)}{3} \ps + (2q-1)(m-1) \ge (2q-1)(m-1)\ge 3(m-1)$ we may assume without
  loss of generality that $p>m-1$ and $p\ge \ps$.
  We then test the first equation in (\ref{0eps}) by $\neps^{p-1}$ and use Young's inequality along with
  (\ref{Deps}), (\ref{D}) and (\ref{F}) to see that for all $t>0$,
  \bas
	\frac{1}{p} \frac{d}{dt} \io \neps^p
	+ (p-1) \kd \io \neps^{p+m-3} |\nabla\neps|^2
	&\le& \frac{1}{p} \frac{d}{dt} \io \neps^p
	+ (p-1) \io \neps^{p-2} D_\eps(\neps) |\nabla\neps|^2 \\
	&=& (p-1) \io \neps^{p-1} F_\eps'(\neps) \nabla\neps\cdot\nabla\ceps \\
	&\le& \frac{(p-1)\kd}{2} \io \neps^{p+m-3} |\nabla\neps|^2
	+ \frac{p-1}{2\kd} \io \neps^{p-m+1} |\nabla\ceps|^2
  \eas
  so that
  \be{8.6}
	\frac{d}{dt} \io \neps^p + C_1 \io \Big|\nabla\neps^\frac{p+m-1}{2}\Big|^2
	\le \frac{p(p-1)}{2\kd} \io \neps^{p-m+1} |\nabla\ceps|^2
	\qquad \mbox{for all } t>0
  \ee
  with $C_1:=\frac{2p(p-1)\kd}{(p+m-1)^2}$.
  Now in order to further estimate the right-hand side herein, we invoke the H\"older inequality to obtain
  \be{8.66}
	\io \neps^{p-m+1} |\nabla\ceps|^2
	\le \bigg\{ \io \neps^{(p-m+1) q'} \bigg\}^\frac{1}{q'} \cdot 
	\bigg\{ \io |\nabla\ceps|^{2q} \bigg\}^\frac{1}{q}
	\qquad \mbox{for all } t>0
  \ee
  with $q':=\frac{q}{q-1}$, where we firstly note that in the case when $(p-m+1)q' \le \ps$, (\ref{8.2}) together
  with the H\"older inquatlity yield $C_2>0$ such that
  \be{8.7}
	\bigg\{ \io \neps^{(p-m+1)q'} \bigg\}^\frac{1}{q'} \le C_2
	\qquad \mbox{for all } t>0.
  \ee
  If, conversely, $(p-m+1)q'>\ps$ then due to our assumption $q\ge 2$ we have
  \bas
	\frac{2(p-m+1)q'}{p+m-1} \le 2q' \le 4<6
  \eas
  and thus $W^{1,2}(\Omega) \hra L^\frac{2(p-m+1)q'}{p+m-1}(\Omega) \hra L^\frac{2\ps}{p+m-1}(\Omega)$,
  whence in particular the number
  \bas
	a:=\frac{3(p+m-1) [(p-m+1)q' -\ps]}{(p-m+1) [3(p+m-1)-\ps] q'}
  \eas
  satisfies $a\in [0,1]$, and accordingly the Gagliardo-Nirenberg inequality provides $C_3>0$ such that
  \bas
	\bigg\{ \io\neps^{(p-m+1)q'} \bigg\}^\frac{1}{q'}
	&=&  \Big\|\neps^\frac{p+m-1}{2}\Big\|_{L^\frac{2(p-m+1)q'}{p+m-1}(\Omega)}^\frac{2(p-m+1)}{p+m-1} \\
	&\le& C_3 \Big\|\nabla\neps^\frac{p+m-1}{2}\Big\|_{L^2(\Omega)}^\frac{6[(p-m+1)q'-\ps]}{[3(p+m-1)-\ps]q'}
	\Big\|\neps^\frac{p+m-1}{2}\Big\|_{L^\frac{2\ps}{p+m-1}(\Omega)}^{\frac{2(p-m+1)}{p+m-1}(1-a)}
	+ C_3 \Big\|\neps^\frac{p+m-1}{2}\Big\|_{L^\frac{2\ps}{p+m-1}(\Omega)}^\frac{2(p-m+1)}{p+m-1}
  \eas
  for all $t>0$. As 
  \bas
	\Big\|\neps^\frac{p+m-1}{2}\Big\|_{L^\frac{2\ps}{p+m-1}(\Omega)}^\frac{2\ps}{p+m-1} = \io \neps^{\ps} 
	\le K \qquad \mbox{for all } t>0
  \eas
  by (\ref{8.2}), together with (\ref{8.7}), (\ref{8.66}) and Young's inequality this shows that regardless of
  the sign of $(p-m+1)q' - \ps$ we can find $C_4>0$ and $C_5>0$ fulfilling
  \bea{8.8}
	\frac{p(p-1)}{2\kd} \io \neps^{p-m+1} |\nabla\ceps|^2
	&\le& C_4 \cdot \bigg\{ 
	\Big\|\nabla\neps^\frac{p+m-1}{2}\Big\|_{L^2(\Omega)}^\frac{6[(p-m+1)q' -\ps]}{[3(p+m-1)-\ps] q'} \, +1 \bigg\}
	\cdot \|\nabla\ceps\|_{L^{2q}(\Omega)}^2 \nn\\
	&\le& 2^{-q'} C_1 \cdot \bigg\{
	\Big\|\nabla\neps^\frac{p+m-1}{2}\Big\|_{L^2(\Omega)}^\frac{6[(p-m+1)q' -\ps]}{[3(p+m-1)-\ps] q'} \, +1 
	\bigg\}^{q'}
	+ C_5 \|\nabla\ceps\|_{L^{2q}(\Omega)}^{2q} \nn\\
	&\le& \frac{C_1}{2} \cdot \bigg\{
	\Big\|\nabla\neps^\frac{p+m-1}{2}\Big\|_{L^2(\Omega)}^\frac{6[(p-m+1)q' -\ps]}{3(p+m-1)-\ps} \, +1 \bigg\}
	+ C_5\|\nabla\ceps\|_{L^{2q}(\Omega)}^{2q}
  \eea
  for all $t>0$, the latter inequality being valid because $(\xi+\eta)^{q'} \le 2^{q'-1} (\xi^{q'}+\eta^{q'})$ for all
  $\xi\ge 0$ and $\eta\ge 0$.\\
  Now our assumption (\ref{8.1}) enters by ensuring that
  \bas
	\frac{6[(p-m+1)q'-\ps]}{3(p+m-1)-\ps} -2
	&=& \frac{6(q'-1)p - 4\ps - 6(m-1)(q'+1)}{3(p+m-1)-\ps} \\
	&=& \frac{6}{[3(p+m-1)-\ps](q-1)} \cdot \Big\{ p- \frac{2(q-1)}{3} \ps - (2q-1)(m-1)\Big\} \\[2mm]
	&\le& 0,
  \eas
  whence another application of Young's inequality yields
  \bas
	\frac{C_1}{2} \cdot \bigg\{
	\Big\|\nabla\neps^\frac{p+m-1}{2}\Big\|_{L^2(\Omega)}^\frac{6[(p-m+1)q' -\ps]}{3(p+m-1)-\ps} \, +1 \bigg\}
	\le \frac{C_1}{2} \io \Big|\nabla\neps^\frac{p+m-1}{2}\Big|^2 + C_1
	\qquad \mbox{for all } t>0.
  \eas
  Together with (\ref{8.8}), this shows that (\ref{8.6}) implies that 
  \be{8.9}
	\frac{d}{dt} \io \neps^p 
	+ \frac{C_1}{2} \io \Big|\nabla\neps^\frac{p+m-1}{2}\Big|^2
	\le C_5 \io |\nabla\ceps|^{2q} + C_1
	\qquad \mbox{for all } t>0,
  \ee
  where a linear absorptive term can be generated again by interpolation in a straightforward manner:
  As according to our restriction $p>\ps$ we know that 
  $W^{1,2}(\Omega)\hra L^\frac{2p}{p+m-1}(\Omega) \hra L^\frac{2\ps}{p+m-1}(\Omega)$, the number
  $b:=\frac{3(p+m-1)(p-\ps)}{[3(p+m-1)-\ps]p}$ satisfies $b\in [0,1]$ and from the Gagliardo-Nirenberg inequality,
  (\ref{8.2}) and Young's inequality we obtain $C_6>0$ and $C_7>0$ such that
  \bas
	\io \neps^p
	&=& \Big\| \neps^\frac{p+m-1}{2}\Big\|_{L^\frac{2p}{p+m-1}(\Omega)}^\frac{2p}{p+m-1} \\
	&\le& C_6 \Big\| \nabla\neps^\frac{p+m-1}{2}\Big\|_{L^2(\Omega)}^\frac{6(p-\ps)}{3(p+m-1)-\ps}
	\Big\|\neps^\frac{p+m-1}{2}\Big\|_{L^\frac{2\ps}{p+m-1}(\Omega)}^{\frac{2p}{p+m-1}(1-b)} 
	+ C_6\Big\|\neps^\frac{p+m-1}{2}\Big\|_{L^\frac{2\ps}{p+m-1}(\Omega)}^\frac{2p}{p+m-1} \\
	&\le& C_7 \Big\|\nabla\neps^\frac{p+m-1}{2}\Big\|_{L^2(\Omega)}^\frac{6(p-\ps)}{3(p+m-1)-\ps} + C_7 \\
	&\le& C_7 \io \Big|\nabla\neps^\frac{p+m-1}{2}\Big|^2 + 2C_7
	\qquad \mbox{for all } t>0,
  \eas
  because $\frac{6(p-\ps)}{3(p+m-1)-\ps} \le \frac{6(p-\ps)}{3p-\ps} \le 2$ by nonnegativity of $m-1$ and $\ps$.
  Therefore, (\ref{8.9}) shows that if we let $y(t):=\io \neps^p(\cdot,t)$, $t\ge 0$, and
  $h(t):=C_5 \io |\nabla\ceps(\cdot,t)|^{2q} + \frac{3}{2} C_1$, $t>0$, than
  \be{8.10}
	y'(t) + \frac{C_1}{4C_7} y(t)
	+ \frac{C_1}{4} \io \Big|\nabla\neps^\frac{p+m-1}{2}\Big|^2 \le h(t)
	\qquad \mbox{for all } t>0,
  \ee
  where in view of our assumption (\ref{8.3}) we have
  \be{8.11}
	\int_t^{t+1} h(s) ds
	\le C_8:=C_5 K + \frac{3}{2} C_1
	\qquad \mbox{for all } t\ge 0.
  \ee
  In view of an elementary lemma on decay in linear first-order ODEs with suitably decaying inhomogeneities
  (see e.g.~\cite[Lemma 3.4]{ssw}), (\ref{8.10}) thus firstly implies that with some $C_9>0$ we have
  $y(t)\le C_9$ for all $t>0$, whereupon (\ref{8.10}) and (\ref{8.11}) secondly entail that
  \bas
	\frac{C_1}{4} \int_t^{t+1} \io \Big|\nabla\neps^\frac{p+m-1}{2}\Big|^2 ds
	\le y(t) + \int_t^{t+1} h(s) ds
	\le C_8+C_9
	\qquad \mbox{for all } t\ge 0,
  \eas
  so that indeed both (\ref{8.4}) and (\ref{8.5}) hold with some conveniently large $C=C(K)>0$.
\qed
\mysection{Uniform $L^p$ bounds on $\neps$ for $p<9(m-1)$ by a first iteration}\label{sect5}
In a first series of applications of Lemma \ref{lem8}, with regard to the regularity assumptions on $\nabla\ceps$
we shall exclusively rely on the corresponding estimate provided by Lemma \ref{lem2} and intend to
repeatedly increase the integrability parameter in (\ref{8.4}) and (\ref{8.5}),
thus keeping the number $q:=2$ in Lemma \ref{lem8} fixed while successively choosing larger values of $\ps$ and $p$.
We shall see that this indeed leads to improved information whenever $m>\frac{10}{9}$, and thereby we
partially re-discover a similar observation that was already made in \cite{taowin_ANIHPC}, with an important 
difference consisting in the fact that unlike in the latter reference, here the achieved bounds are global in time.
\begin{lem}\label{lem6}
  Let $m>\frac{10}{9}$. Then for all $p\in [1,9(m-1))$ there exists $C(p)>0$ such that for all $\eps\in (0,1)$,
  \be{6.1}
	\io \neps^p(\cdot,t) \le C(p)
	\qquad \mbox{for all } t\ge 0
  \ee
  and
  \be{6.01}
	\int_t^{t+1} \io \Big|\nabla\neps^\frac{p+m-1}{2}\Big|^2 \le C(p)
	\qquad \mbox{for all } t\ge 0.
  \ee
\end{lem}
\proof
  We define $(p_k)_{k\in\N_0} \subset\R$ by letting $p_0:=1$ and 
  \be{6.11}
	p_{k+1}:=\frac{2}{3} p_k + 3(m-1)
	\qquad \mbox{for } k\ge 0.
  \ee
  It can the readily be verified that due to our assumption $m>\frac{10}{9}$ the sequence $(p_k)_{k\in\N_0}$ is 
  strictly increasing with $p_k\nearrow 9(m-1)$ as $k\to\infty$, so that by means of an interpolation argument
  it is clear that we only need to prove (\ref{6.1}) and (\ref{6.01}) for $p=p_k$ and each $k\in\N_0$.
  To this end, we note that the case $k=0$ can be covered by combining Lemma \ref{lem2} with (\ref{mass}), so that
  in view of an inductive reasoning we are left with the verification of the property that whenever $k\in\N_0$ is such 
  that
  \be{6.2}
	\io \neps^{p_k}(\cdot,t) \le C_1(k)
	\quad \mbox{and} \quad
	\int_t^{t+1} \io \Big|\nabla\neps^\frac{p_k+m-1}{2}\Big|^2 \le C_1(k)
	\qquad \mbox{for all $t\ge 0$ and each } \eps\in (0,1)
  \ee
  with some $C_1(k)>0$, we can find $C_2(k)>0$ satisfying
  \be{6.3}
	\io \neps^{p_{k+1}}(\cdot,t) \le C_2(k)
	\quad \mbox{and} \quad
	\int_t^{t+1} \io \Big|\nabla\neps^\frac{p_{k+1}+m-1}{2}\Big|^2 \le C_2(k)
	\qquad \mbox{for all $t\ge 0$ and any } \eps\in (0,1).
  \ee
  To achieve this, we observe that according to the first inequality in (\ref{6.2}) and (\ref{2.2}), the requirements
  (\ref{8.2}) and (\ref{8.4}) from Lemma \ref{lem8} are fulfilled for $\ps:=p_k$ and $q:=2$.
  In light of (\ref{6.11}), both inequalities in (\ref{6.3}) therefore result from an application of Lemma \ref{lem8}
  to $p:=p_{k+1}$.
\qed
\mysection{Improving estimates for $\nabla\ceps$ via maximal Sobolev regularity}\label{sect6}
We next plan to apply Lemma \ref{lem8} by using the outcome of Lemma \ref{lem6} as a starting point with respect
to the regularity assumptions on $\neps$, but with regard to the hypothesis (\ref{8.3}) no longer going back to
Lemma \ref{lem2} but rather using suitably improved integrability information on $\nabla\ceps$. 
Within a range of $m$ which is smaller than that in Lemma \ref{lem6} but yet larger than the interval 
$(\frac{9}{8},\infty)$ we shall finally focus on, such further properties can indeed be gained under the assumptions
provided by the result of Lemma \ref{lem6} by means of the key Lemma \ref{lem5} below which in turn relies
on the following statement on time-independent bounds for $\ueps$ in appropriate Lebesgue spaces.
\begin{lem}\label{lem3}
  Let $m>\frac{215}{192}$. Then there exists
  $\delta_1(m)>0$ such that for all $p>1$ fulfilling
  $p>9(m-1)-\delta_1(m)$ and $K>0$ one can find $C(p,K)>0$ with the property that if for some $\eps\in (0,1)$ we have
  \be{3.1}
	\io \neps^p(\cdot,t) \le K
	\qquad \mbox{for all } t\ge 0,
  \ee
  then
  \be{3.2}
	\io |\ueps(\cdot,t)|^\frac{2(5p+3m-3)}{3} \le C(p,K)
	\qquad \mbox{for all } t \ge 0.
  \ee
\end{lem}
\proof
  We let
  \bas
	\rho(p):=20p^2-(33-12m)p-18(m-1), \qquad p\in\R.
  \eas
  Then our assumption $m>\frac{215}{192}$ precisely warrants that
  \bas
	\rho(9(m-1))
	= 1620(m-1)^2 -9\cdot(33-12m)(m-1) - 18(m-1) 
	= 9(m-1)(192m-215) 
	>0,
  \eas
  while since $\frac{215}{192}>\frac{131}{124}$ we moreover have
  \bas
	\rho'(p)
	=40 p-33+12m 
	\ge 360 (m-1) - 33+12m
	= 3\cdot (124m-131) >0
	\qquad \mbox{for all } p>9(m-1).
  \eas
  We can therefore pick $\delta_1=\delta_1(m)>0$ such that
  \bas
	\rho(p)>0
	\qquad \mbox{for all } p>9(m-1)-\delta_1(m),
  \eas
  and given $p>1$ such that $p>9(m-1)-\delta_1(m)$ we thus obtain that $q:=\frac{2(5p+3m-3)}{3}$ satisfies $q>1$ and
  \bas
	3(3-2p) \cdot \Big(q-\frac{3p}{3-2p}\Big)
	= - \rho(p)<0
  \eas
  and hence 		
  \be{3.3}
	\frac{3}{2}\Big(\frac{1}{p}-\frac{1}{q}\Big)<1.
  \ee
  Now assuming (\ref{3.1}) for some $\eps\in (0,1)$ and $K>0$, on the basis of a variation-of-constants representation
  of $\ueps$ we can estimate
  \bea{3.4}
	\|\ueps(\cdot,t)\|_{L^q(\Omega)}
	\le \|e^{-tA} u_0\|_{L^q(\Omega)}
	+ \int_0^t \Big\| e^{-(t-s)A} \proj[ \neps(\cdot,s)\nabla \phi] \Big\|_{L^q(\Omega)} ds,
	\qquad t>0,
  \eea
  and recall known regularization properties of the Dirichlet Stokes semigroup $(e^{-tA})_{t\ge 0}$
  (\cite[p.201]{giga1986}) to find
  $C_1>0, C_2>0$ and $\lambda>0$ such that
  \be{3.5}
	\|e^{-tA} u_0\|_{L^q(\Omega)}
	\le C_1 \|u_0\|_{L^q(\Omega)}
	\qquad \mbox{for all } t>0
  \ee
  and
  \bea{3.6}
	\int_0^t \Big\| e^{-(t-s)A} \proj[ \neps(\cdot,s)\nabla \phi] \Big\|_{L^q(\Omega)} ds
	\le C_2 \int_0^t \Big(1+(t-s)^{-\frac{3}{2}(\frac{1}{p}-\frac{1}{q})} \Big) e^{-\lambda(t-s)} 
	\Big\| \proj[\neps(\cdot,s)\nabla\phi] \Big\|_{L^p(\Omega)} ds
  \eea
  for all $t>0$.
  Here by boundedness of $\nabla\phi$ on $\Omega$ and the continuity of the Helmholtz projection when acting as an 
  operator in $L^p(\Omega;\R^3)$ (\cite{fujiwara_morimoto}), we see that with some $C_3>0$ we have
  \bas
	\Big\|\proj[\neps(\cdot,s)\nabla\phi]\Big\|_{L^p(\Omega)}
	\le C_3\|\neps(\cdot,s)\|_{L^p(\Omega)} \le C_3 K^\frac{1}{p}
	\qquad \mbox{for all } s>0
  \eas
  according to (\ref{3.1}). Therefore, (\ref{3.6}) entails that
  \bas
	\int_0^t \Big\| e^{-(t-s)A} \proj[ \neps(\cdot,s)\nabla \phi] \Big\|_{L^q(\Omega)} ds
	&\le& C_2 C_3 K^\frac{1}{p} \int_0^t \Big(1+(t-s)^{-\frac{3}{2}(\frac{1}{p}-\frac{1}{q})}\Big) e^{-\lambda(t-s)} 
	ds \nn\\
	&\le& C_2 C_3 C_4 K^\frac{1}{p}
	\qquad \mbox{for all } t>0
  \eas 
  with $C_4:=\int_0^\infty (1+\sigma^{-\frac{3}{2}(\frac{1}{p}-\frac{1}{q})}) e^{-\lambda\sigma} d\sigma$ being
  finite thanks to (\ref{3.3}). 
  When combined with (\ref{3.5}) and (\ref{3.4}), in view of our choice of $q$ this establishes (\ref{3.2}).
\qed
As a second preliminary for Lemma \ref{lem8},
let us note how a pair of hypotheses in the flavor of (\ref{bound})
influences space-time integrability of $\neps$ by means of straighforward interpolation.
\begin{lem}\label{lem4}
  Let $m>1$. Then for all $p\ge 1$ and any $K>0$ there exists $C(p,K)>0$ such that if for some $\eps\in (0,1)$ we have
  \be{4.1}
	\io \neps^p(\cdot,t) \le K
	\qquad \mbox{for all } t\ge 0
  \ee
  and
  \be{4.2}
	\int_t^{t+1} \io \Big|\nabla\neps^\frac{p+m-1}{2} \Big|^2 \le K
	\qquad \mbox{for all } t\ge 0,
  \ee
  then
  \be{4.3}
	\int_t^{t+1} \io \neps^\frac{5p+3m-3}{3}
	\le C(p,K)
	\qquad \mbox{for all } t\ge 0.
  \ee
\end{lem}
\proof
  Using that $p>0$ and $m\ge 1$ imply that
  \bas
	\frac{2p}{p+m-1} \le \frac{2(5p+3m-3)}{3(p+m-1)} \le 6
  \eas
  and hence $W^{1,2}(\Omega) \hra L^\frac{2(5p+3m-3)}{3(p+m-1)}(\Omega) \hra L^\frac{2p}{p+m-1}(\Omega)$,
  from the Gagliardo-Nirenberg inequality we obtain $C_1>0$ such that
  \bas
	\io \neps^\frac{5p+3m-3}{3}
	&=& \Big\| \neps^\frac{p+m-1}{2}\Big\|_{L^\frac{2(5p+3m-3)}{3(p+m-1)}(\Omega)}^\frac{2(5p+3m-3)}{3(p+m-1)} \\
	&\le& C_1 \Big\|\nabla\neps^\frac{p+m-1}{2}\Big\|_{L^2(\Omega)}^2
	\Big\| \neps^\frac{p+m-1}{2}\Big\|_{L^\frac{2p}{p+m-1}(\Omega)}^\frac{4p}{3(p+m-1)}
	+ C_1\Big\| \neps^\frac{p+m-1}{2}\Big\|_{L^\frac{2p}{p+m-1}(\Omega)}^\frac{2(5p+3m-3)}{3(p+m-1)}
	\qquad \mbox{for all } t>0.
  \eas
  Noting that $\|\neps^\frac{p+m-1}{2}\|_{L^\frac{2p}{p+m-1}(\Omega)}^\frac{2p}{p+m-1} \le K$ 
  for all $t>0$ by (\ref{4.1}),
  on integrating in time we thus infer that
  \bas
	\int_t^{t+1} \io \neps^\frac{5p+3m-3}{3}
	&\le& C_1 K^\frac{2}{3} \int_t^{t+1} \io \Big|\nabla\neps^\frac{p+m-1}{2}\Big|^2
	+ C_1 K^\frac{5p+3m-3}{3p} \\
	&\le& C_1 K^\frac{5}{3} 
	+ C_1 K^\frac{5p+3m-3}{3p}
  \eas
  for all $t\ge 0$.
\qed
We can now proceed to the main result of this section which, on the basis of a maximal regularity
property of scalar parabolic equations, asserts that bounds of the flavor in (\ref{bound})
entail an estimate for $\nabla\ceps$ in a spatio-temporal $L^{2q}$ space with some positive $q$ 
which indeed satisfies $q>2$ if $p\ge 1$ is suitably large.
\begin{lem}\label{lem5}
  Let $m>\frac{215}{192}$, and let $\delta_1(m)>0$ be as in Lemma \ref{lem3}. Then for all $p>9(m-1)-\delta_1(m)$
  and each $K>0$ one can find $C(p,K)>0$ with the property that if for some $\eps\in (0,1)$,
  \be{5.1}
	\io \neps^p(\cdot,t) \le K
	\qquad \mbox{for all }  \ge >0
  \ee
  and
  \be{5.2}
	\int_t^{t+1} \io \Big|\nabla\neps^\frac{p+m-1}{2} \Big|^2 \le K
	\qquad \mbox{for all } t\ge 0,
  \ee
  then
  \be{5.3}
	\int_t^{t+1} \io |\nabla\ceps|^\frac{2(5p+3m-3)}{3} \le C(p,K)
	\qquad \mbox{for all } t\ge 0.
  \ee
\end{lem}
\proof
  We abbreviate $q:=\frac{5p+3m-3}{3}$ and apply a standard result on maximal Sobolev regularity in scalar parabolic
  equations (\cite{giga_sohr}) to find $C_1>0$ with the property that whenever $t_\star\in\R$,
  $z\in C^{2,1}(\bom\times [t_\star,t_\star+2])$ and $f\in C^0(\bom\times [\ts,\ts+2])$ are such that
  \bas
	\left\{ \begin{array}{ll}
	z_t=\Delta z + f(x,t), 
	\qquad & x\in\Omega, \ t\in (\ts,\ts+2), \\[1mm]
	\frac{\partial z}{\partial\nu}=0,
	\qquad & x\in\pO, \ t\in (\ts,\ts+2), \\[1mm]
	z(x,\ts)=0,
	\qquad & x\in\Omega, 
	\end{array} \right.
  \eas
  then
  \be{5.4}
	\int_{\ts}^{\ts+2} \|z(\cdot,t)\|_{W^{2,q}(\Omega)}^q dt
	\le C_1 \int_{\ts}^{\ts+2} \|f(\cdot,t)\|_{L^q(\Omega)}^q dt.
  \ee
  Furthermore, let us fix $C_2>0$ and $C_3>0$ such that in accordance with a well-known regularization feature of
  the Neumann heat semigroup (\cite{win_JDE}) and the Gagliardo-Nirenberg inequality we have
  \be{5.444}
	\|\nabla e^{t\Delta} \varphi\|_{L^{2q}(\Omega)}
	\le C_2\|\varphi\|_{W^{1,2q}(\Omega)}
	\qquad \mbox{for all $\varphi\in W^{1,2q}(\Omega)$ and any } t>0
  \ee
  as well as
  \be{5.44}
	\|\nabla\varphi\|_{L^{2q}(\Omega)}^{2q}
	\le C_3 \|\varphi\|_{W^{2,q}(\Omega)}^q \|\varphi\|_{L^\infty(\Omega)}^q
	\qquad \mbox{for all } \varphi\in W^{2,q}(\Omega),
  \ee
  where in establishing the latter we note that $W^{2,q}(\Omega) \hra W^{1,2q}(\Omega) \hra L^\infty(\Omega)$ due to
  the fact that $q\ge \frac{5}{3}>\frac{3}{2}$.\\
  As a final preparation, let us observe that according to Lemma \ref{lem4} and Lemma \ref{lem3}, our assumptions 
  (\ref{5.1}) and (\ref{5.2}) ensure that we can choose $C_4(K)>0$ and $C_5(K)>0$ such that
  \be{5.45}
	\int_t^{t+2} \|\neps(\cdot,s)\|_{L^q(\Omega)}^q ds \le C_4(K)
	\qquad \mbox{for all } t\ge 0
  \ee
  and
  \be{5.46}
	\|\ueps(\cdot,t)\|_{L^{2q}(\Omega)} \le C_5(K)
	\qquad \mbox{for all } t\ge 0,
  \ee
  the latter conclusion relying on our hypothesis on $p$.\abs
  In order to make appropriate use of these preliminaries in the present context, we pick a nondecreasing
  $\zeta_0\in C^\infty(\R)$ such that $\zeta_0\equiv 0$ in $(-\infty,-1]$ and $\zeta_0\equiv 1$ in $[1,\infty)$,
  and for fixed $t_0\ge 0$ we let $\zeta(t)\equiv \zeta^{(t_0)}(t):=\zeta_0(t-t_0)$, $t\ge 0$, and
  \bas
	z(\cdot,t):=\zeta(t)\cdot \Big\{ \ceps(\cdot,t)- e^{t\Delta} c_0\Big\},
	\qquad t \ge (t_0-1)_+.
  \eas
  Then by (\ref{0eps}) and the identity $\partial_t e^{t\Delta} c_0=\Delta e^{t\Delta} c_0$,
  \bea{5.5}
	z_t
	&=& \zeta(t)\cdot \Big\{ \Delta\ceps - F_\eps(\neps)\ceps - \ueps\cdot\nabla\ceps \Big\}
	-\zeta(t)\Delta e^{t\Delta} c_0
	+ \zeta'(t) \cdot \Big\{ \ceps - e^{t\Delta} c_0\Big\} \nn\\
	&=& \Delta z - \zeta(t) F_\eps(\neps)\ceps - \ueps\cdot\nabla z \nn\\
	& & - \zeta(t) \ueps\cdot \nabla e^{t\Delta} c_0
	+ \zeta'(t) \ceps - \zeta'(t) e^{t\Delta} c_0
	\qquad \mbox{in } \Omega\times ((t_0-1)_+,\infty),
  \eea
  and clearly 
  \be{5.6}
	\frac{\partial z}{\partial\nu}=0
	\qquad \mbox{on } \pO\times ((t_0-1)_+,\infty).
  \ee
  Moreover, at the respective initial time we have
  \be{5.7}
	z\Big(\cdot,(t_0-1)_+\Big) \equiv 0
	\qquad \mbox{in } \Omega,
  \ee
  because if $t_0\ge 1$ then $\zeta(t_0-1)=0$ and hence 
  \bas
	z\Big(\cdot,(t_0-1)_+\Big)
	= z(\cdot,t_0-1)
	= \zeta(t_0-1) \cdot \Big\{\ceps(\cdot,t_0-1) - e^{(t_0-1)\Delta} c_0\Big\} =0
	\qquad \mbox{in } \Omega,
  \eas
  whereas if $t_0\in [0,1)$ then
  \bas
	z\Big(\cdot,(t_0-1=_+\Big)
	= z(\cdot,0)
	= \zeta(0) \cdot \Big\{ \ceps(\cdot,0)-c_0\Big\} =0
	\qquad \mbox{in } \Omega
  \eas
  by (\ref{0eps}).\\
  As a consequence of (\ref{5.5})-(\ref{5.7}), we may now invoke (\ref{5.4}) which along with (\ref{5.44}) and
  (\ref{cinfty}) shows that abbreviating $\ts:=(t_0-1)_+$ we have
  \bea{5.8}
	\int_{\ts}^{\ts+2} \|\nabla z(\cdot,t)\|_{L^{2q}(\Omega)}^{2q} dt
	&\le& C_3 \|c_0\|_{L^\infty(\Omega)}^q \int_{\ts}^{\ts+2} \|z(\cdot,t)\|_{W^{2,q}(\Omega)}^q dt \nn\\
	&\le& 5^q C_1 C_3 \|c_0\|_{L^\infty(\Omega)}^q
	\int_{\ts}^{\ts+2} \bigg\{ \|\zeta(t) F_\eps(\neps) \ceps\|_{L^q(\Omega)}^q
	+ \|\ueps\cdot\nabla z\|_{L^q(\Omega)}^q \nn\\[2mm]
	& & \hspace*{40mm}
	+ \|\zeta(t) \ueps\cdot\nabla e^{t\Delta} c_0\|_{L^q(\Omega)}^q
	+ \|\zeta'(t)\ceps\|_{L^q(\Omega)}^q \nn\\[2mm]
	& & \hspace*{40mm}
	+ \|\zeta'(t) e^{t\Delta} c_0\|_{L^q(\Omega)}^q \bigg\} \, dt.
  \eea
  Here we use that by (\ref{F}) we have $0\le F_\eps(s)\le s$ for all $s\ge 0$ and that $0\le \zeta(t)\le 1$
  for all $t\in\R$ to see, again by means of (\ref{cinfty}), that  
  \bea{5.9}
	\int_{\ts}^{\ts+2} \|\zeta(t) F_\eps(\neps)\ceps\|_{L^q(\Omega)}^q dt
	&\le& \|c_0\|_{L^\infty(\Omega)}^q
	\int_{\ts}^{\ts+2} \|\neps(\cdot,t)\|_{L^q(\Omega)}^q dt \nn\\
	&\le& C_4(K) \|c_0\|_{L^\infty(\Omega)}^q
  \eea
  according to (\ref{5.45}), while the Cauchy-Schwarz ineaulity together with (\ref{5.46}) and (\ref{5.444}) shows that
  \bea{5.10}
	\int_{\ts}^{\ts+2} \|\zeta(t)\ueps\cdot\nabla e^{t\Delta} c_0\|_{L^q(\Omega)}^q dt
	&\le& \int_{\ts}^{\ts+2} \|\ueps(\cdot,t)\|_{L^{2q}(\Omega)}^q
	\|\nabla e^{t\Delta} c_0\|_{L^{2q}(\Omega)}^q dt \nn\\
	&\le& 2C_2^q C_5^q(K) \|c_0\|_{W^{1,2q}(\Omega)}^q.
  \eea
  Next, by (\ref{cinfty}) and the contractivity of the semigroup $(e^{t\Delta})_{t\ge 0}$ on $L^q(\Omega)$,
  writing $C_6:=\|\zeta_0'\|_{L^\infty(\R)}$ we obtain
  \be{5.11}
	\int_{\ts}^{\ts+2} \|\zeta'(t)\ceps\|_{L^q(\Omega)}^q dt
	\le 2C_6^q |\Omega|\cdot \|c_0\|_{L^\infty(\Omega)}^q
  \ee
  and
  \be{5.12}
	\int_{\ts}^{\ts+2} \|\zeta'(t) e^{t\Delta} c_0\|_{L^q(\Omega)}^q dt
	\le 2C_6^q \|c_0\|_{L^q(\Omega)}^q,
  \ee
  so that it remains to estimate the corresponding integral associated with the second summand in brackets on the right
  of (\ref{5.8}).
  For this purpose, after employing the Cauchy-Schwarz inequality we additionally make use of Young's inequality to see,
  again by means of (\ref{5.46}), that
  \bas
	C_1 C_3\|c_0\|_{L^\infty(\Omega)}^q
	\int_{\ts}^{\ts+2} \|\ueps\cdot\nabla z\|_{L^q(\Omega)}^q dt 
	&\le&
	C_1 C_3\|c_0\|_{L^\infty(\Omega)}^q
	\int_{\ts}^{\ts+2} \|\ueps\|_{L^{2q}(\Omega)}^q \|\nabla z\|_{L^{2q}(\Omega)}^q dt \\
	&\le& \frac{1}{2} \int_{\ts}^{\ts+2} \|\nabla z\|_{L^{2q}(\Omega)}^{2q} dt
	+ \frac{C_1^2 C_3^2 \|c_0\|_{L^\infty(\Omega)}^{2q}}{2} 
	\int_{\ts}^{\ts+2} \|\ueps\|_{L^{2q}(\Omega)}^{2q} dt \\
	&\le& \frac{1}{2} \int_{\ts}^{\ts+2} \|\nabla z\|_{L^{2q}(\Omega)}^{2q} dt
	+ C_1^2 C_3^2 C_5^{2q}(K) \|c_0\|_{L^\infty(\Omega)}^{2q}.
  \eas
  In conjunction with (\ref{5.9})-(\ref{5.12}), this shows that (\ref{5.8}) leads to the inequality
  \bas
	\frac{1}{2} \int_{\ts}^{\ts+2} \|\nabla z(\cdot,t)\|_{L^{2q}(\Omega)}^{2q} dt
	&\le& C_7(K):= C_4(K) \|c_0\|_{L^\infty(\Omega)}^q
	+ C_1^2 C_3^2 C_5^{2q}(K) \|c_0\|_{L^\infty(\Omega)}^{2q} \\
	& & \hspace*{16mm}
	+ 2C_2^q C_5^q(K) \|c_0\|_{W^{1,2q}(\Omega)}^q
	+ 2C_6^q \|c_0\|_{L^\infty(\Omega)}^q
	+ 2C_6^q \|c_0\|_{L^q(\Omega)}^q,
  \eas
  so that since $(t_0,t_0+1) \subset \Big((t_0-1)_+ , (t_0-1)_++2\Big)$ and thus $\zeta\equiv 1$ in
  $(t_0,t_0+1)$, in particular we infer that
  \bas
	\int_{t_0}^{t_0+1} \Big\| \nabla\ceps(\cdot,t) - \nabla e^{t\Delta} c_0\Big\|_{L^{2q}(\Omega)}^{2q} dt
	\le 2C_7(K)
	\qquad \mbox{for all } t_0\ge 0.
  \eas
  Once more recalling (\ref{5.444}) and using that $(\xi+\eta)^{2q} \le 2^{2q-1} (\xi^{2q}+\eta^{2q})$ for all
  $\xi\ge 0$ and $\eta\ge 0$, we therefore obtain that
  \bas
	\int_{t_0}^{t_0+1} \|\nabla\ceps(\cdot,t)\|_{L^{2q}(\Omega)}^{2q} dt
	&\le& 2^{2q-1} \int_{t_0}^{t_0+1} \Big\|\nabla\ceps(\cdot,t)-\nabla e^{t\Delta} c_0\Big\|_{L^{2q}(\Omega)}^{2q} dt
	+ 2^{2q-1} \int_{t_0}^{t_0+1} \|\nabla e^{t\Delta} c_0\|_{L^{2q}(\Omega)}^{2q} dt \\
	&\le& 2^{2q} C_7(K)
	+ 2^{2q-1} C_2^{2q} \|c_0\|_{W^{1,2q}(\Omega)}^{2q}
	\qquad \mbox{for all } t_\ge 0,
  \eas
  which in view of our definition of $q$ precisely yields (\ref{5.3}).
\qed
\mysection{Arbitrary $L^p$ bounds for $\neps$ by a second iteration}\label{sect7}
Now in light of Lemma \ref{lem5}, our general regularity statement from Lemma \ref{lem8} can readily developed
to the following basis for a second iterative reasoning.
\begin{lem}\label{lem9}
  Let $m>\frac{215}{192}$ and $\ps>9(m-1)-\delta_1(m)$ with $\delta_1(m)>0$ taken from Lemma \ref{lem3}.
  Then for all $p>1$ fulfilling
  \be{9.1}
	p\le \frac{10\ps^2 + (36m-42)\ps + (m-1)(18m-27)}{9},
  \ee
  and any choice of $K>0$ one can pick $C(p,K)>0$ such that if for some $\eps\in (0,1)$ we have
  \be{9.2}
	\io \neps^{\ps}(\cdot,t) \le K
	\qquad \mbox{for all } t\ge 0
  \ee
  and
  \be{9.3}
	\int_t^{t+1} \io \Big|\nabla\neps^\frac{\ps+m-1}{2}\Big|^2 \le K
	\qquad \mbox{for all } t\ge 0,
  \ee
  then
  \be{9.4}
	\io \neps^p(\cdot,t) \le C(p,K)
	\qquad \mbox{for all } t\ge 0
  \ee
  and
  \be{9.5}
	\int_t^{t+1} \io \Big|\nabla\neps^\frac{p+m-1}{2}\Big|^2 \le C(p,K)
	\qquad \mbox{for all } t\ge 0.
  \ee
\end{lem}
\proof
  Since $\ps>9(m-1)-\delta_1(m)$, we may invoke Lemma \ref{lem5} to see that writing $q:=\frac{5\ps+3m-3}{3}$ we can find
  $C_1(K)>0$ such that
  \be{9.6}
	\int_t^{t+1} \io |\nabla\ceps|^{2q} \le C_1(K)
	\qquad \mbox{for all } t\ge 0.
  \ee
  Now observing that in our situation the right-hand side of (\ref{8.1}) can be rewritten according to
  \bas
	\frac{2(q-1)}{3} \ps + (2q-1)(m-1)
	&=& \frac{2\cdot \frac{5\ps+3m-6}{3}}{3} \cdot \ps
	+ \frac{10\ps+6m-9}{3} \cdot (m-1) \\
	&=& \frac{10\ps^2 + (36m-42)\ps + (m-1)(18m-27)}{9},
  \eas
  given any $p>1$ fulfilling (\ref{9.1}) we may apply Lemma \ref{lem8} to infer that due to (\ref{9.2}) and (\ref{9.6})
  both inequalities in (\ref{9.4}) and (\ref{9.5}) hold if we fix $C(p,K)>0$ suitably large.
\qed
With regard to the question how far the above lemma through its condition (\ref{9.1}) indeed allows for an improvement 
in knowledge, let us briefly prove the following elementary observations which highlight the role
of the restriction $m>\frac{9}{8}$ made in Theorem \ref{theo99}.
\begin{lem}\label{lem10}
  For $m>1$, let
  \be{10.1}
	\psi(p):=\frac{10p^2 + (36m-42)p + (m-1)(18m-27)}{9},
	\qquad p\in\R.
  \ee
  Then 
  \be{10.2}
	\psi\Big(9(m-1)\Big) >9(m-1)
	\quad \mbox{if and only if} \quad m>\frac{9}{8},
  \ee
  and there exist $\delta_2(m)>0$ and $\Gamma>1$ such that
  \be{10.3}
	\psi(p) \ge\Gamma p
	\qquad \mbox{for all } p>9(m-1)-\delta_2(m).
  \ee
\end{lem}
\proof
  Computing
  \bas
	\frac{\psi(9(m-1)) - 9(m-1)}{m-1}
	&=& \frac{810(m-1)^2 + 9(36m-42)(m-1) + (m-1)(18m-27)}{9(m-1)} - 9 \\[2mm]
	&=& 16(8m-9) \\[2mm]
	&>& 0,
  \eas
  we directly obtain (\ref{10.2}).
  To verify (\ref{10.3}), we let 
  \bas
	\tpsi(p):=\frac{\psi(p)}{p}
	\qquad \mbox{for } p>0,
  \eas
  so that since (\ref{10.2}) asserts that $C_1:=\tpsi(9(m-1))-1$ is positive, by continuity we can pick  
  $\delta_2=\delta_2(m)>0$ such that $9(m-1)-\delta_2>0$ and
  \be{10.4}
	\tpsi(p) \ge \Gamma:=1+\frac{C_1}{2}
	\qquad \mbox{for all } p \in \Big( 9(m-1)-\delta_2, 9(m-1)\Big].
  \ee
  As
  \be{10.5}
	\tpsi'(p) = \frac{10}{9} - \frac{(m-1)(2m-3)}{p^2}
	\qquad \mbox{for all } p>0,
  \ee
  it thus immediately follows that if $m\le \frac{3}{2}$ then $\tpsi'\ge\frac{10}{9}>0$ throughout $(0,\infty)$.
  If $m>\frac{3}{2}$, then for $p\ge 9(m-1)$ we can use (\ref{10.5}) to estimate
  \bas
	\tpsi'(p) \ge \frac{10}{9} - \frac{(m-1)(2m-3)}{81(m-1)^2}
	= \frac{88 m - 87}{81(m-1)} >0,
  \eas
  because $m>1$.
  In both cases, we thus obtain that $\tpsi'>0$ on $[9(m-1),\infty)$ and hence $\tpsi\ge \Gamma$ on
  $(9(m-1)-\delta_2,\infty)$ by (\ref{10.4}).
\qed
We are thereby prepared for our second recursive argument, with its outcome being as follows.
\begin{lem}\label{lem11}
  Let $m>\frac{9}{8}$. Then for all $p>1$ there exists $C(p)>0$ such that
  \be{11.1}
	\io \neps^p(\cdot,t) \le C(p)
	\qquad \mbox{for all } t\ge 0
  \ee
  and
  \be{11.01}
	\int_t^{t+1} \io \Big|\nabla\neps^\frac{p+m-1}{2}\Big|^2 \le C(p)
	\qquad \mbox{for all } t\ge 0.
  \ee
\end{lem}
\proof
  As $m>\frac{9}{8}>\frac{215}{192}$, taking $\delta_1(m)>0$ and $\delta_2(m)>0$ as given by Lemma \ref{lem3} and
  Lemma \ref{lem10}, respectively, we may pick $p_0\in (1,9(m-1))$ such that
  \be{11.2}
	p_0>9(m-1) - \min\{\delta_1(m),\delta_2(m)\},
  \ee
  and thereupon recursively define
  \be{11.21}
	p_k:=\psi(p_{k-1}),
	\qquad k\in\N=\{1,2,3,...\},
  \ee
  with $\psi:\R\to\R$ taken from Lemma \ref{lem10}. 
  Then since $p_0>9(m-1)-\delta_2(m)$ by (\ref{11.2}), according to (\ref{10.2}) an inductive argument shows that
  \be{11.22}
	p_k \ge \Gamma^k p_0
	\qquad \mbox{for all } k\in\N
  \ee
  with $\Gamma>1$ as provided by Lemma \ref{lem10}, whence in particular $p_k\to\infty$ as $k\to\infty$.
  Now due to the boundedness of $\Omega$, in order to verify the lemma it is sufficient to show that for all $k\ge 0$
  there exists $C_1(k)>0$ such that for all $\eps\in (0,1)$,
  \be{11.3}
	\io \neps^p(\cdot,t) \le C_1(k)
	\quad \mbox{and} \quad
	\int_t^{t+1} \io \Big|\nabla\neps^\frac{p_k+m-1}{2}\Big|^2 \le C_1(k)
	\qquad \mbox{for all } t\ge 0,
  \ee
  which will again result from an iterative reasoning:
  Namely, for $k=0$ the claimed inequality is a direct consequence of Lemma \ref{lem6}, because 
  $m>\frac{9}{8}>\frac{10}{9}$ and $p_0\in (1,9(m-1))$.
  If (\ref{11.3}) holds for some $k_0\ge 0$ and some $C_1(k_0)>0$, however, then since (\ref{11.22}) and (\ref{11.2})
  warrant that $p_k\ge p_0>9(m-1)-\delta_1(m)$, and again since $m>\frac{215}{192}$, Lemma \ref{lem9} provides
  $C_2>0$ such that
  \bas
	\io \neps^p(\cdot,t) \le C_2
	\quad \mbox{and} \quad
	\int_t^{t+1} \io \Big|\nabla\neps^\frac{p+m-1}{2}\Big|^2 \le C_2
	\qquad \mbox{for all } t\ge 0
  \eas
  with
  \bas
	p:=\frac{10 p_{k_0}^2 + (36m-42) p_{k_0} + (m-1)(18m-27)}{9}.
  \eas
  As thus $p=\psi(p_{k_0})=p_{k_0+1}$ by (\ref{11.21}), this asserts (\ref{11.3}) also for $k=k_0+1$
  and thereby completes the proof.
\qed
\mysection{Further regularity properties}\label{sect8}
With Lemma \ref{lem11} at hand, further regularity properties can now be obtained by essentially straightforward 
arguments: We firstly recall Lemma \ref{lem3} and a standard regularization feature of the heat semigroup
to obtain the following.
\begin{lem}\label{lem21}
  Let $p>1$. Then there exists $C(p)>0$ such that whenever $\eps\in (0,1)$,
  \be{21.2}
	\io |\nabla\ceps(\cdot,t)|^p \le C(p)
	\qquad \mbox{for all } t\ge 0
  \ee
  and
  \be{21.1}
	\io |\ueps(\cdot,t)|^p \le C(p)
	\qquad \mbox{for all } t\ge 0.
  \ee
\end{lem}
\proof
  In view of Lemma \ref{lem11}, (\ref{21.1}) is an evident consequence of Lemma \ref{lem3}.
  Thereafter, (\ref{21.2}) can be derived from (\ref{21.1}) and again Lemma \ref{lem11} by well-known
  results on gradient regularity in semilinear heat equations (\cite{horstmann_win}).
\qed
By means of a Moser iteration, the latter together with Lemma \ref{lem11} entails an $\eps$-independent
$L^\infty$ bound for $\neps$.
\begin{lem}\label{lem22}
  There exists $C>0$ such that for arbitrary $\eps\in (0,1)$ we have
  \be{22.1}
	\|\neps(\cdot,t)\|_{L^\infty(\Omega)} \le C
	\qquad \mbox{for all } t\ge 0.
  \ee
\end{lem}
\proof
  In view of Lemma \ref{lem21} and Lemma \ref{lem11} when applied to suitably large $p>1$, this directly
  follows from a Moser-type iterative procedure (see \cite[Lemma A.1]{taowin_JDE} for a version precisely covering
  the present case).
\qed
Again by means of maximal Sobolev regularity properties combined with an appropriate embedding result, 
the estimates collected above imply H\"older bounds for $\ceps,\ueps$ and $\nabla\ceps$.
This will be achieved in Lemma \ref{lem42} on the basis of the following lemma in which any influence
of the respective initial data is faded out.
\begin{lem}\label{lem41}
  There exist $\theta\in (0,1)$ and $C>0$ such that for all $\eps\in (0,1)$,
  \be{41.1}
	\|\ceps-\hc\|_{C^{1+\theta,\theta}(\bom\times [t,t+1])} \le C
	\qquad \mbox{for all } t\ge 0
  \ee
  and
  \be{41.2}
	\|\ueps-\hu\|_{C^{1+\theta,\theta}(\bom\times [t,t+1])} \le C
	\qquad \mbox{for all } t\ge 0,
  \ee
  where
  \be{41.3}
	\hc(\cdot,t):=e^{t\Delta} c_0
	\quad \mbox{and} \quad
	\hu(\cdot,t):=e^{-tA} u_0
	\qquad \mbox{for } t\ge 0.
  \ee
\end{lem}
\proof
  Since $\hc_t=\Delta\hc$, it follows from (\ref{0eps}) that
  \bas
	\partial_t (\ceps-\hc) = \Delta(\ceps-\hc) - F_\eps(\neps) \ceps - \ueps\cdot\nabla\ceps,
	\qquad x\in\Omega, \ t>0,
  \eas
  where given $p>1$ we may invoke Lemma \ref{lem21}, Lemma \ref{lem22} and (\ref{cinfty}) and recall (\ref{F})
  to find $C_1>0$ fulfilling
  \bas
	\int_t^{t+2} \io \Big|-F_\eps(\neps)\ceps - \ueps\cdot\nabla\ceps \Big|^p \le C_1
	\qquad \mbox{for all $t\ge 0$ and } \eps\in (0,1).
  \eas
  Therefore, by means of maximal Sobolev regularity estimates along with an appropriate time localization in the style of 
  the argument from Lemma \ref{lem5}, we infer the existence of $C_2>0$ such that
  \bas
	\int_t^{t+2} \Big\{ \|\ceps(\cdot,s)-\hc(\cdot,s)\|_{W^{2,p}(\Omega)}^p
	+ \|\partial_t(\ceps(\cdot,s)-\hc(\cdot,s))\|_{L^p(\Omega)}^p \Big\} ds \le C_2
	\qquad \mbox{for all $t\ge 0$ and } \eps\in (0,1).
  \eas
  In view of a known embedding property (\cite{amann}), an application thereof to suitably large $p>1$ establishes
  (\ref{41.1}).\\
  Likewise, using that
  \bas
	\partial_t(\ueps-\hu)=-A(\ueps-\hu) + \proj[\neps\nabla\phi],
	\qquad x\in\Omega, \ t>0,
  \eas
  and that herein for $p>1$ we can use the boundedness of $\proj$ on $L^p(\Omega;\R^3)$ (\cite{fujiwara_morimoto})
  together with Lemma \ref{lem22} to find $C_3>0$ such that
  \bas
	\int_t^{t+2} \io \Big|\proj[\neps(\cdot,s)\nabla\phi]\Big|^p \le C_3
	\qquad \mbox{for all $t\ge 0$ and } \eps\in (0,1),
  \eas
  we obtain (\ref{41.2}) from corresponding maximal Sobolev regularity estimates for the Stokes evolution
  equation (\cite{giga_sohr}).
\qed
Indeed, the latter inter alia implies the following H\"older estimates, which with regard to the gradient bound
in (\ref{42.3}) must remain local in time due to possibly lacking appropriate regularity and compatibility properties
of $c_0$.
\begin{lem}\label{lem42}
  There exists $\theta\in (0,1)$ with the property that one can find $C>0$ such that for all $\eps\in (0,1)$,
  \be{42.1}
	\|\ceps\|_{C^\theta(\bom\times [t,t+1])} \le C
	\qquad \mbox{for all } t\ge 0
  \ee
  and
  \be{42.2}
	\|\ueps\|_{C^\theta(\bom\times [t,t+1])} \le C
	\qquad \mbox{for all } t\ge 0,
  \ee
  and that for all $\tau>0$ it is possible to choose $C(\tau)>0$ fulfilling
  \be{42.3}
	\|\nabla\ceps\|_{C^\theta(\bom\times [t,t+1])} \le C(\tau)
	\qquad \mbox{for all } t\ge \tau
  \ee
  whenever $\eps\in (0,1)$.
\end{lem}
\proof
  We take $\hc$ and $\hu$ from (\ref{41.3}) and note that since $c_0\in W^{1,\infty}(\Omega) \hra
  \bigcap_{\theta\in (0,1)} C^\theta(\bom)$ and
  $u_0\in D(A^\alpha) \hra \bigcap_{\theta\in (0,2\alpha-\frac{3}{2})} C^\theta(\bom)$
  (\cite{giga1981_theother}, \cite{henry}), known smoothing properties of the heat equation 
  and the Stokes evolution system
  ensure that there exist $\theta_1\in (0,1), \theta_2\in (0,1), C_1>0$ and $C_2>0$ such that
  \bas
	\|\hc\|_{C^{\theta_1}(\bom\times [t,t+1])} \le C_1
	\qquad \mbox{for all } t\ge 0
  \eas
  and
  \bas
	\|\hu\|_{C^{\theta_2}(\bom\times [t,t+1])} \le C_2
	\qquad \mbox{for all } t\ge 0,
  \eas
  and that for all $\tau>0$ we can find $C_3(\tau)>0$ such that
  \bas
	\|\nabla\hc\|_{C^1(\bom\times [t,t+1])} \le C_3(\tau)
	\qquad \mbox{for all } t\ge\tau.
  \eas
  Therefore, (\ref{42.1})-(\ref{42.3}) result from Lemma \ref{lem41}.
\qed
For strongly degenerate cell diffusion present when e.g.~$D(s)=s^{m-1}$, $s\ge 0$, with
large values of $m$, we do not know whether $\neps$ enjoys equicontinuity properties in the classical pointwise sense,
which may indeed suffer from a possible dominance of the transport terms in the first equation of (\ref{0eps}) at small
densities.
In order to nevertheless provide some compactness and equicontinuity properties of this solution component,
let us finally derive two statements on time regularity of $\neps$ in a straightforward manner.
\begin{lem}\label{lem23}
  Let $T>0$. Then there exists $C(T)>0$ such that for all $\eps\in (0,1)$,
  \be{23.1}
	\int_0^T \Big\|\partial_t \neps^m(\cdot,t)\Big\|_{(W_0^{1,\infty}(\Omega))^\star} dt \le C(T)
  \ee
  and
  \be{23.11}
	\|n_{\eps t}(\cdot,t)\|_{(W_0^{2,2}(\Omega))^\star} \le C(T)
	\qquad \mbox{for all } t\in (0,T). 
  \ee
\end{lem}
\proof
  We fix $t\in (0,T)$ and 
  $\zeta\in C_0^\infty(\Omega)$ such that $\|\zeta\|_{W^{1,\infty}(\Omega)} \le 1$, and then obtain from
  the first equation in (\ref{0eps}) by straightforward manipulations that writing
  $C_1:=\sup_{\eps\in (0,1)} \|\neps\|_{L^\infty(\Omega\times (0,\infty))}$ and
  $C_2:=\|D\|_{L^\infty((0,c_1))} +2$, according to (\ref{D}) we have
  \bas
	& & \hspace*{-18mm}
	\bigg| \frac{1}{m} \io \partial_t \neps^m(\cdot,t) \zeta \bigg| \\
	&=& \bigg| \io \neps^{m-1} \nabla\cdot \Big\{ D_\eps(\neps) \nabla\neps
	- \neps F_\eps'(\neps) \nabla\ceps
	- \neps\ueps \Big\} \zeta \bigg| \\
	&=& \bigg| - (m-1) \io \neps^{m-2} D_\eps(\neps) |\nabla\neps|^2 \zeta
	- \io \neps^{m-1} D_\eps(\neps) \nabla\neps\cdot\nabla\zeta \\
	& & \hspace*{5mm}
	+ (m-1) \io \neps^{m-1} F_\eps'(\neps) (\nabla\neps\cdot\nabla\ceps) \zeta
	+ \io \neps^m F_\eps'(\neps) \nabla\ceps\cdot\nabla\zeta 
	+ \frac{1}{m} \io \neps^m \ueps\cdot\nabla\zeta \bigg| \\
	&\le& (m-1) C_2 \io \neps^{m-2} |\nabla\neps|^2
	+ (m-1)C_2 \io \neps^{m-1} |\nabla\neps| \\
	& & + (m-1) \io \neps^{m-1} |\nabla\neps| \cdot |\nabla\ceps|
	+ \io \neps^m |\nabla\ceps|
	+ \frac{1}{m} \io \neps^m |\ueps| \\
	&\le& (m-1) C_2 \io \neps^{m-2} |\nabla\neps|^2
	+ \frac{m-1}{2} C_2 \io \neps^{m-2} |\nabla\neps|^2
	+ \frac{m-1}{2} C_2 \io \neps^m \\
	& & + \frac{m-1}{2} \io \neps^{m-2} |\nabla\neps|^2
	+ \frac{m-1}{2} \io \neps^m |\nabla\ceps|^2 \\
	& & + \io \neps^m |\nabla\ceps|
	+ \frac{1}{m} \io \neps^m |\ueps| \\
	&\le& (m-1)\Big( \frac{3C_2}{2} + \frac{1}{2}\Big) \io \neps^{m-2} |\nabla\neps|^2
	+ \frac{(m-1) C_1^m C_2 |\Omega|}{2} + \frac{(m-1) C_1^m}{2} \io |\nabla\ceps|^2 \\
	& & + C_1^m \io |\nabla\ceps|
	+ \frac{C_1^m}{m} \io |\ueps|
	\qquad \mbox{for all } \eps\in (0,1).
  \eas
  In view of the estimates provided by Lemma \ref{lem2} and Lemma \ref{lem21}, (\ref{23.1}) therefore readily results 
  upon integration.\\
  The inequality in (\ref{23.11}) can similarly be derived from Lemma \ref{lem21} and Lemma \ref{lem22}.
\qed
\mysection{Existence of a global bounded weak solution}\label{sect9}
In the sequel, we shall refer to the following natural concept of weak solvability in (\ref{0}), (\ref{0i}), (\ref{0b}):
\begin{defi}\label{defi_weak}
  Let 
  \bea{reg_w1}
	& & n \in L^1_{loc}(\bom \times [0,\infty)), \nn\\
	& & c \in L^\infty_{loc}(\bom \times [0,\infty)) \cap L^1_{loc}([0,\infty);W^{1,1}(\Omega)) 
	\qquad \mbox{and} \nn\\
	& & u \in L^1_{loc}([0,\infty);W^{1,1}(\Omega;\R^3)),
  \eea
  be such that $n\ge 0$ and $c\ge 0$ in $\Omega\times (0,T)$ and
  \bea{reg_w2}
	D_0(n), \ n|\nabla c| \ \mbox{and} \ n|u|
	\ \mbox{ belong to } L^1_{loc}(\bom\times [0,\infty)),
  \eea
  where $D_0(s):=\int_0^s D(\sigma)d\sigma$ for $s\ge 0$.
  Then $(n,c,u)$ will be called a {\em global weak solution} of (\ref{0}), (\ref{0i}), (\ref{0b}) if
  $\nabla \cdot u=0$ in the distributional sense, if
  \be{w1}
	- \int_0^\infty \io n\varphi_t - \io n_0 \varphi(\cdot,0)
	= \int_0^\infty \io D_0(n) \Delta \varphi
	+ \int_0^\infty \io n \nabla c\cdot \nabla \varphi
	+ \int_0^\infty \io nu\cdot \nabla \varphi
  \ee
  for all $\varphi\in C_0^\infty(\bom\times [0,\infty))$ fulfilling $\frac{\partial\varphi}{\partial\nu}=0$ on
  $\partial\Omega\times (0,\infty)$, if
  \be{w2}
	- \int_0^\infty \io c\varphi_t - \io c_0 \varphi(\cdot,0)
	= - \int_0^\infty \io \nabla c \cdot \nabla \varphi
	- \int_0^\infty \io nc\varphi
	+ \int_0^\infty \io cu\cdot \nabla \varphi
  \ee
  for all $\varphi\in C_0^\infty(\bom\times [0,\infty))$, and if moreover
  \be{w3}
	- \int_0^\infty \io u\cdot \varphi_t - \io u_0\cdot \varphi(\cdot,0)
	= - \int_0^\infty \io \nabla u \cdot \nabla \varphi
	+ \int_0^\infty \io n\nabla \phi \cdot \varphi
  \ee
  for all $\varphi \in C_0^\infty(\Omega\times [0,\infty);\R^3)$ 
  such that $\nabla \varphi\equiv 0$ in $\Omega\times (0,\infty)$.
\end{defi}
In this context, a series of standard extraction procedures on the basis of our estimates collected above
indeed yields global solvability.
\begin{lem}\label{lem24}
  Let $m>\frac{9}{8}$. Then there exist $(\eps_j)_{j\in\N} \subset (0,1)$, a null set $N\subset (0,\infty)$
  and a triple $(n,c,u)$ of functions $n:\Omega\times (0,\infty)\to [0,\infty)$,
  $c:\Omega\times (0,\infty)\to [0,\infty)$ and $u:\Omega\times (0,\infty)\to \R^3$ such that
  $\eps_j\searrow 0$ as $j\to\infty$ and
  \begin{eqnarray}
	& & \neps(\cdot,t) \to n(\cdot,t)
	\qquad \mbox{a.e.~in $\Omega$ for all } t\in (0,\infty)\setminus N, 
	\label{24.1} \\
	& & \neps \wsto n
	\qquad \mbox{in } L^\infty(\Omega\times (0,\infty)),
	\label{24.2} \\
	& & \neps \to n
	\qquad \mbox{in } C^0_{loc}([0,\infty);(W_0^{2,2}(\Omega))^\star),
	\label{24.3} \\
	& & \ceps \to c
	\qquad \mbox{in } C^0_{loc}(\bom\times [0,\infty)),
	\label{24.4} \\
	& & \ceps\wsto c
	\qquad \mbox{in } L^\infty((0,\infty);W^{1,p}(\Omega))
	\qquad \mbox{for all } p\in (1,\infty),
	\label{24.6} \\
	& & \nabla\ceps\to\nabla c
	\qquad \mbox{in } C^0_{loc}(\bom\times [0,\infty)),
	\label{24.5} \\
	& & \ueps\to u
	\qquad \mbox{in } C^0_{loc}(\bom\times [0,\infty)),
	\label{24.7} \\
	& & \ueps\wsto u
	\qquad \mbox{in } L^\infty(\Omega\times (0,\infty))
	\qquad \mbox{and}
	\label{24.77} \\
	& & \nabla\ueps\to\nabla u
	\qquad \mbox{in } L^2_{loc}(\bom\times [0,\infty))
	\label{24.8}
  \eea
  as $\eps=\eps_j\searrow 0$.
  Moreover, $(n,c,u)$ forms a global weak solution of (\ref{0}), (\ref{0i}), (\ref{0b})
  in the sense of Definition \ref{defi_weak}, and we have
  \be{24.9}
	\io n(\cdot,t) = \io n_0
	\qquad \mbox{for all } t\in (0,\infty)\setminus N.
  \ee
\end{lem}
\proof
  Since Lemma \ref{lem2}, Lemma \ref{lem22} and Lemma \ref{lem23} guarantee that $(\neps^m)_{\eps\in (0,1)}$
  is bounded in $L^2_{loc}([0,\infty;W^{1,2}(\Omega))$ and that
  $(\partial_t \neps^m)_{\eps\in (0,1)}$ is bounded in $L^2_{loc}([0,\infty);(W_0^{3,2}(\Omega))^\star)$
  due to the continuity of the embedding $W_0^{3,2}(\Omega) \hra W_0^{1,\infty}(\Omega)$, 
  an Aubin-Lions lemma (\cite{temam}) yields $(\eps_j)_{j\in\N}\subset (0,1)$ such that
  $\eps_j\searrow 0$ as $j\to\infty$ and that $\neps^m\to n^m$ holds a.e.~in $\Omega\times (0,\infty)$ as
  $\eps=\eps_j\searrow 0$ with some nonnegative function $n$ defined on $\Omega\times (0,\infty)$, whence
  using the Fubini-Tonelli theorem we readily obtain (\ref{24.1}).
  In view of Lemma \ref{lem22}, Lemma \ref{lem2} and (\ref{23.11}), on further extraction we may also achieve
  (\ref{24.2}) and (\ref{24.3}), whereas the bounds provided by Lemma \ref{lem2}, Lemma \ref{lem21} and
  Lemma \ref{lem42} ensure that we can moreover easily achieve (\ref{24.4})-(\ref{24.8}) upon two applications
  of the Arzel\`a-Ascoli theorem.\\
  The regularity properties in (\ref{reg_w1}) and (\ref{reg_w2}) as well as the claimed solenoidality of $u$
  are evident from (\ref{24.1})-(\ref{24.8}),
  while the verification of (\ref{w1}), (\ref{w2}) and (\ref{w3}) is thereafter straightforward.
\qed
\mysection{Large time behavior}\label{sect10}
\subsection{Basic decay information}
Next addressing the large time asymptotics of our solutions,
as in several previous studies on 
qualitative behavior in related chemotaxis-fluid
systems with signal absorption (\cite{win_ARMA}, \cite{lankeit_M3AS}, \cite{win_TRAN}, \cite{win_CVPDE})
we shall rely on the following elementary information indicating a certain decay of the quantities $nc$ and $\nabla c$.
Here and throughout the sequel, without further mentioning we shall assume that $m>\frac{9}{8}$ and that
$(n,c,u)$ denotes the global weak solution constructed in Lemma \ref{lem24}.
\begin{lem}\label{lem25}
  There exist $\eps_\star\in (0,1)$ and $C>0$ such that 
  \be{25.1}
	\int_0^\infty \io \neps\ceps \le C
	\qquad \mbox{for all } \eps\in (0,\eps_\star)
  \ee
  and
  \be{25.2}
	\int_0^\infty \io |\nabla\ceps|^2 \le C
	\qquad \mbox{for all } \eps\in (0,\eps_\star).
  \ee
\end{lem}
\proof
  Using Lemma \ref{lem22}, we can fix $C_1>0$ such that $\neps\le C_1$ in $\Omega\times (0,\infty)$ for all
  $\eps\in (0,1)$, and let $\eps_\star\in (0,1)$ be small enough such that $\frac{1}{\eps_\star}\ge C_1$.
  Then (\ref{def_F}) implies that $F_\eps(\neps)\equiv \neps$ throughout $\Omega\times (0,\infty)$ whenever
  $\eps\in (0,\eps_\star)$, whence integrating the second equation in (\ref{0eps}) we obtain
  \bas
	\io \ceps(\cdot,t) + \int_0^t \io \neps\ceps = \io c_0
	\qquad \mbox{for all $\eps\in (0,\eps_\star)$ and each } t>0,
  \eas
  from which (\ref{25.1}) follows.
  Moreover, testing the same equation by $\ceps$ and recalling (\ref{F}) yields
  \bas
	\frac{1}{2} \io \ceps^2(\cdot,t) + \int_0^t \io |\nabla\ceps|^2 
	= \frac{1}{2} \io c_0^2 - \int_0^t \io F_\eps(\neps)\ceps
	\le \frac{1}{2} \io c_0^2
	\qquad \mbox{for all $\eps\in (0,1)$ and } t>0
  \eas
  and thereby verifies (\ref{25.2}).
\qed
\subsection{Decay of $c$}
A first application of Lemma \ref{lem25} shows that thanks to the uniform H\"older estimates from Lemma \ref{lem42}
the second solution component indeed decays in the sense claimed in Theorem \ref{theo99}.
\begin{lem}\label{lem32}
  We have
  \be{32.1}
	c(\cdot,t) \to 0
	\quad \mbox{in } W^{1,\infty}(\Omega)
	\qquad \mbox{as } t\to\infty.
  \ee
\end{lem}
\proof
  Following a variant of an approach pursued in \cite{win_ARMA}, we first use (\ref{24.9}) and the Poincar\'e
  inequality to see that for all $\eps\in (0,1)$,
  \bas
	\onz \cdot \io \ceps
	&=& \io \neps \overline{\ceps} \\
	&=& \io \neps\ceps - \io \neps (\ceps-\overline{\ceps}) \\
	&\le& \io \neps\ceps
	+ \sqrt{|\Omega|}C_1 \bigg\{ \io (\ceps-\overline{\ceps})^2 \bigg\}^\frac{1}{2} \\
	&\le& \io \neps\ceps
	+ C_2 \bigg\{ \io |\nabla\ceps|^2 \bigg\}^\frac{1}{2}
	\qquad \mbox{for all } t>0
  \eas
  with $C_1:=\sup_{\eps\in (0,1)} \|\neps\|_{L^\infty(\Omega\times (0,\infty))}<\infty$ by Lemma \ref{lem22},
  and with some $C_2>0$.
  Thus, by (\ref{cinfty}),
  \bas
	\onz^2 \cdot \bigg\{ \io \ceps\bigg\}^2
	&\le& 2\bigg\{ \io \neps\ceps\bigg\}^2 + 2C_2^2 \io |\nabla\ceps|^2 \\
	&\le& 2C_1 \|c_0\|_{L^\infty(\Omega)} \io \neps\ceps + 2C_2^2 \io |\nabla\ceps|^2
	\qquad \mbox{for all } t>0,
  \eas
  so that according to Lemma \ref{lem25} we infer that with some $\eps_\star\in (0,1)$ and $C_3>0$ we have
  \bas
	\int_0^\infty \|\ceps(\cdot,t)\|_{L^1(\Omega)}^2 dt  \le C_3
	\qquad \mbox{for all } \eps\in (0,\eps_\star)
  \eas
  and hence
  \bas
	\int_0^\infty \|c(\cdot,t)\|_{L^1(\Omega)}^2 dt \le C_3
  \eas
  thanks to Lemma \ref{lem24} and Fatou's lemma.
  Since the spatio-temporal H\"older continuity property expressed by (\ref{42.1}) warrants that
  $0\le t \mapsto \|c(\cdot,t)\|_{L^1(\Omega)}$ is uniformly continuous, through a standard argument
  this entails that necessarily
  \be{32.2}
	c(\cdot,t)\to 0
	\quad \mbox{in } L^1(\Omega)
	\qquad \mbox{as } t\to\infty.
  \ee
  Since Lemma \ref{lem42} moreover guarantees that with some $\theta\in (0,1)$ and $C_4>0$ we have
  \be{32.3}
	\|c(\cdot,t)\|_{C^{1+\theta}(\bom)} \le C_4
	\qquad \mbox{for all } t>1,
  \ee
  a straightforward reasoning based on interpolation and the compactness of the first among the continuous embeddings
  $C^{1+\theta}(\bom) \hra W^{1,\infty}(\Omega) \hra L^1(\Omega)$ shows that (\ref{32.2}) 
  and (\ref{32.3}) entail (\ref{32.1}):
  In fact, given $\eta>0$ we may employ an Ehrling-type lemma to pick $C_5>0$ fulfilling
  \be{32.4}
	\|\varphi\|_{W^{1,\infty}(\Omega)} 
	\le \frac{\eta}{2C_4} \|\varphi\|_{C^{1+\theta}(\bom)} + C_5\|\varphi\|_{L^1(\Omega)}
	\qquad \mbox{for all } \varphi\in C^{1+\theta}(\bom),
  \ee
  and then use (\ref{32.2}) to choose $t_0>1$ satisfying
  \bas
	\|c(\cdot,t)\|_{L^1(\Omega)} \le \frac{\eta}{2C_5}
	\qquad \mbox{for all } t>t_0.
  \eas
  Then by (\ref{32.4}) and (\ref{32.3}),
  \bas
	\|c(\cdot,t)\|_{L^1(\Omega)} \le \frac{\eta}{2C_4} \cdot C_4 + C_5 \cdot \frac{\eta}{2C_5}=\eta
	\qquad \mbox{for all } t>t_0,
  \eas
  as desired.
\qed
\subsection{Stabilization of $n$}
Next concerned with the large time behavior of $n$, in order to circumvent obstacles stemming from 
possibly strong degeneracies of diffusion when $m$ is large, we rely on another quasi-energy structure
in deriving the following result which can be viewed as asserting a certain short-time conservation 
of smallness of the quantity $\io (\neps-\onz)^2$, and which, remarkably, beyond the above properties and in particular
(\ref{25.2}) does not explicitly require the presence of any diffusion mechanism in the first equation in (\ref{0eps}).
\begin{lem}\label{lem31}
  There exists $C>0$ such that for each $\eps\in (0,1)$ and any choice of $\ts\ge 0$ we have
  \bea{31.1}
	\io \Big(\neps(\cdot,t)-\onz\Big)^2
	&\le& C\cdot \bigg\{ \io \Big(\neps(\cdot,\ts)-\onz\Big)^2
	+ \io |\nabla\ceps(\cdot,\ts)|^2
	+ \sup_{s\in (\ts,\ts+1)} \|\ceps(\cdot,s)\|_{L^2(\Omega)}^2 \bigg\} \nn\\[2mm]
	& & \hspace*{60mm}
	\qquad \mbox{for all } t\in (\ts,\ts+1).
  \eea
\end{lem}
\proof
  We start by multiplying the first equation in (\ref{0eps}) by $\neps-\onz$ to obtain
  \bea{31.2}
	\frac{1}{2} \frac{d}{dt} \io (\neps-\onz)^2
	&=& - \io D_\eps(\neps) |\nabla\neps|^2 
	+ \io \neps F_\eps'(\neps) \nabla\neps\cdot\nabla\ceps \nn\\
	&\le& \io \neps F_\eps'(\neps) \nabla\neps\cdot\nabla\ceps
	\qquad \mbox{for all } t>0.
  \eea
  Here in order to appropriately estimate the right-hand side, we introduce
  \bas
	G_\eps(s):=\int_0^s \sigma F_\eps'(\sigma) d\sigma, \qquad s\ge 0,
  \eas
  and once more integrate by parts to rewrite
  \bea{31.3}
	\io \neps F_\eps'(\neps) \nabla\neps\cdot\nabla\ceps
	&=& \io \nabla G_\eps(\neps)\cdot\nabla\ceps \nn\\
	&=& - \io G_\eps(\neps) \Delta\ceps \nn\\
	&=& - \io \Big(G_\eps(\neps)-G_\eps(\onz)\Big)\cdot \Delta\ceps
	\qquad \mbox{for all } t>0,
  \eea
  because $\io \Delta\ceps(\cdot,t)=0$ for all $t>0$.
  Now since we know from Lemma \ref{lem22} that with some $C_1>0$ we have
  \be{31.4}
	\neps \le C_1
	\quad \mbox{in } \Omega\times (0,\infty)
	\qquad \mbox{for all } \eps\in (0,1),
  \ee
  and since $0\le G_\eps'(s) \le s$ thanks to (\ref{F}), by the mean value theorem we can estimate
  \bas
	\Big|G_\eps(\neps(x,t))-G_\eps(\onz)\Big|
	&\le& \|G_\eps'\|_{L^\infty((0,C_1))} |\neps(x,t)-\onz| \\
	&\le& C_1 |\neps(x,t)-\onz|
	\qquad \mbox{for all $x\in\Omega, t>0$ and } \eps\in (0,1).
  \eas
  By means of Young's inequality, (\ref{31.3}) therefore implies that
  \bas
	\io \neps F_\eps'(\neps) \nabla\neps\cdot\nabla\ceps
	&\le& C_1 \io |\neps-\onz| \cdot |\Delta\ceps| \nn\\
	&\le& \frac{1}{2} \io (\neps-\onz)^2
	+ \frac{C_1^2}{2} \io |\Delta\ceps|^2 
	\qquad \mbox{for all } t>0,
  \eas
  and that in view of (\ref{31.2}) we thus have
  \be{31.5}
	\frac{d}{dt} \io (\neps-\onz)^2
	\le \io (\neps-\onz)^2
	+ C_1^2 \io |\Delta\ceps|^2
	\qquad \mbox{for all } t>0.
  \ee
  Here an adequate compensation of the rightmost integral can be achieved by using the second equation in (\ref{0eps}),
  which when tested against $-\Delta\ceps$ yields
  \bea{31.6}
	\frac{1}{2} \frac{d}{dt} \io |\nabla\ceps|^2
	+ \io |\Delta\ceps|^2
	&=& \io F_\eps(\neps) \ceps\Delta\ceps
	+ \io (\ueps\cdot\nabla\ceps) \Delta\ceps \nn\\
	&\le& \frac{1}{4} \io |\Delta\ceps|^2
	+ \io \neps^2 \ceps^2 \nn\\
	& & +\frac{1}{4} \io |\Delta\ceps|^2
	+ \io |\ueps\cdot\nabla\ceps|^2 \nn\\
	&\le& \frac{1}{2} \io |\Delta\ceps|^2
	+ C_1^2 \io \ceps^2 + C_2^2 \io |\nabla \ceps|^2
	\qquad \mbox{for all } t>0,
  \eea
  where in accordance with Lemma \ref{lem42} we have chosen $C_2>0$ large enough fulfilling
  $|\ueps| \le C_2$ in $\Omega\times (0,\infty)$ for all $\eps\in (0,1)$.\\
  In combination, (\ref{31.5}) and (\ref{31.6}) now show that
  \bas
	\frac{d}{dt} \bigg\{ \io (\neps-\onz)^2 + C_1^2 \io |\nabla\ceps|^2 \bigg\}	
	\le \io (\neps-\onz)^2 + 2C_1^4 \io \ceps^2
	+ 2C_1^2 C_2^2 \io |\nabla\ceps|^2
	\qquad \mbox{for all } t>0,
  \eas
  implying that $y(t):=\io (\neps(\cdot,t)-\onz)^2 + C_1^2 \io |\nabla\ceps(\cdot,t)|^2$, $t\ge 0$, satisfies
  \bas
	y'(t) \le C_3 y(t) + C_4 \io \ceps^2
	\qquad \mbox{for all } t>0
  \eas
  with $C_3:=\max\{1, 2C_2^2\}$ and $C_4:=2C_1^4$. By an ODE comparison, this entails that  
  \bas	
	y(t) &\le& e^{C_3(t-\ts)} y(\ts)
	+ C_4 \int_{\ts}^t e^{C_3(t-s)} \cdot \bigg\{ \io \ceps^2(\cdot,s) \bigg\} ds \\
	&\le& e^{C_3} y(\ts)
	+ \frac{C_4 e^{C_3}}{C_3} \cdot \sup_{s\in (\ts,\ts+1)} \io \ceps^2(\cdot,s)
	\qquad \mbox{for all } t\in (\ts,\ts+1)
  \eas
  and thereby establishes (\ref{31.1}).
\qed
By means of another $L^p$ testing procedure applied to the first equation in (\ref{0eps}), again relying on
the estimate (\ref{25.2}) from Lemma \ref{lem25}, the latter implies stabilization of $n$ toward its average,
at least when yet considered in $l^2(\Omega)$ and outside a null set of times.
\begin{lem}\label{lem35}
  Let $N\subset (0,\infty)$ be as provided by Lemma \ref{lem24}.
  Then 
  \be{35.1}
	n(\cdot,t) \to \onz
	\quad \mbox{in } L^2(\Omega)
	\qquad \mbox{as } (0,\infty)\setminus N \ni t\to\infty.
  \ee
\end{lem}
\proof
  We first invoke Lemma \ref{lem31} to find $C_1>0$ such that for any $\ts\ge 0$ and $\eps\in (0,1)$ we have
  \bas
	\io \Big(\neps(\cdot,t)-\onz\Big)^2
	&\le& C_1 \cdot \bigg\{ \io \Big(\neps(\cdot,\ts)-\onz\Big)^2 
	+ \io |\nabla\ceps(\cdot,\ts)|^2
	+ \sup_{s\in (\ts,\ts+1)} \|\ceps(\cdot,s)\|_{L^2(\Omega)}^2 \bigg\}\\[2mm]
	& & \hspace*{60mm}
	\qquad \mbox{for all } t\in (\ts,\ts+1).
  \eas
  Here since $\ceps\to c$ in $C^0_{loc}(\bom\times [0,\infty))$ and 
  $\nabla\ceps\to\nabla c$ in $C^0_{loc}(\bom\times [1,\infty))$as $\eps=\eps_j\searrow 0$ according to Lemma \ref{lem24},
  and since $(\neps-\onz)_{\eps\in (0,1)}$ is bounded in $L^\infty(\Omega\times (0,\infty))$ by Lemma \ref{lem22},
  on the basis of (\ref{24.1}) and the dominated convergence theorem we may let $\eps=\eps_j\searrow 0$ to obtain that
  \bea{35.3}
	\io \Big(n(\cdot,t)-\onz\Big)^2
	&\le& C_1 \cdot \bigg\{ \io \Big(n(\cdot,\ts)-\onz\Big)^2
	+ \io  |\nabla c(\cdot,\ts)|^2
	+ \sup_{s\in (\ts,\ts+1)} \|c(\cdot,s)\|_{L^2(\Omega)}^2 \bigg\} \nn\\[2mm]
	& & \hspace*{20mm}
	\mbox{for all $\ts\in (1,\infty)\setminus N$ and any } t\in (\ts,\ts+1)\setminus N.
  \eea
  In order to prepare an appropriate control of the right-hand side herein, we fix some $\gamma\ge 1$ satisfying 
  $\gamma\ge m-1$ and use $\neps^{2\gamma-m}$ as a test function in the first equation from (\ref{0eps}) to see that
  for all $\eps \in (0,1)$,
  \bas
	\frac{1}{2\gamma-m+1} \frac{d}{dt} \io \neps^{2\gamma-m+1}
	&+& (2\gamma-m) \io \neps^{2\gamma-m-1} D_\eps(\neps) |\nabla \neps|^2 \\
	&=& (2\gamma-m) \io \neps^{2\gamma-m} F_\eps'(\neps) \nabla\neps\cdot\nabla\ceps
	\qquad \mbox{for all } t>0,
  \eas
  which in light of (\ref{Deps}), (\ref{D}), (\ref{F}) and Young's inequality implies that
  \bas
	& & \hspace*{-20mm}
	\frac{1}{2\gamma-m+1} \io \neps^{2\gamma-m+1}(\cdot,t)
	+ \frac{(2\gamma-m)\kd}{2} \int_0^t \io \neps^{2\gamma-2} |\nabla\neps|^2 \nn\\
	&\le& \frac{1}{2\gamma-m+1} \io n_0^{2\gamma-m+1}
	- \frac{(2\gamma-m)\kd}{2} \int_0^t \io \neps^{2\gamma-2} |\nabla\neps|^2  \nn\\
	& & + (2\gamma-m) \int_0^t \io \neps^{2\gamma-m} |\nabla\neps| \cdot |\nabla\ceps| \nn\\
	&\le& \frac{1}{2\gamma-m+1} \io n_0^{2\gamma-m+1}
	+ \frac{2\gamma-m}{2\kd} \int_0^t \io \neps^{2\gamma-2m+2} |\nabla\ceps|^2 \nn\\
	&\le& \frac{1}{2\gamma-m+1} \io n_0^{2\gamma-m+1}
	+ \frac{2\gamma-m}{2\kd} \|\neps\|_{L^\infty(\Omega\times (0,\infty))}^{2\gamma-2m+2}
	\int_0^\infty \io |\nabla\ceps|^2
	\qquad \mbox{for all } t>0,
  \eas
  because $2\gamma-2m+2\ge 0$.
  Due to the boundedness properties asserted by Lemma \ref{lem22} and Lemma \ref{lem25}, we therefore conclude
  that there exist $\eps_\star\in (0,1)$ and $C_2>0$ such that
  \bas
	\int_0^\infty \io \Big|\nabla\neps^\gamma\Big|^2 \le C_2
	\qquad \mbox{for all } \eps\in (0,\eps_\star),
  \eas
  and that hence according to the Poincar\'e inequality we can find $C_3>0$ fulfilling
  \bas
	\int_0^\infty \Big\| \neps^\gamma(\cdot,t)-\mu_\eps^\gamma(t)\Big\|_{L^2(\Omega)}^2 dt \le C_3
	\qquad \mbox{for all } \eps\in (0,\eps_\star),
  \eas
  where we have set $\mu_\eps(t):=\Big\{ \frac{1}{|\Omega|} \io \neps^\gamma(\cdot,t)\Big\}^\frac{1}{\gamma}$ for 
  $t>0$ and $\eps\in (0,1)$.
  Again using (\ref{24.1}) along with the dominated convergence theorem, from this we readily infer on invoking 
  Fatou's lemma on the time interval $(0,\infty)$ that
  \be{35.4}
	\int_0^\infty \Big\| n^\gamma(\cdot,t)-\mu^\gamma(t)\Big\|_{L^2(\Omega)}^2 dt \le C_3
  \ee
  is valid with $\mu(t):=\Big\{ \frac{1}{|\Omega|} \io n^\gamma(\cdot,t)\Big\}^\frac{1}{\gamma}$,
  $t\in (0,\infty)\setminus N$, the latter satisfying $\mu(t)\ge \onz$ for all $t\in (0,\infty)\setminus N$
  due to the fact that by (\ref{24.9}) and the H\"older inequality we can estimate
  \bas
	\onz |\Omega| = \io n(\cdot,t) \le \bigg\{ \io n^\gamma(\cdot,t)\bigg\}^\frac{1}{\gamma}
	\cdot |\Omega|^{1-\frac{1}{\gamma}}
	\qquad \mbox{for } t\in (0,\infty)\setminus N.
  \eas
  As $|\xi^\gamma-\eta^\gamma| \ge \eta^{\gamma-1} \cdot |\xi-\eta|$ for all $\xi\ge 0$ and $\eta\ge 0$, this implies that
  \bas
	\Big|n^\gamma(\cdot,t)-\mu^\gamma(t) \Big|^2
	&\ge& \mu^{2\gamma-2}(t) \cdot |n(\cdot,t)-\mu(t)|^2 \\
	&\ge& \onz^{2\gamma-2} \cdot |n(\cdot,t)-\mu(t)|^2
	\quad \mbox{a.e.~in $\Omega$ for all } t\in (0,\infty)\setminus N,
  \eas
  so that from (\ref{35.4}) we obtain that
  \be{35.5}
	\int_0^\infty \|n(\cdot,t)-\mu(t)\|_{L^2(\Omega)}^2 dt \le C_4
  \ee
  with $C_4:=C_3 \cdot \onz^{2-2\gamma}>0$.\\
  Now to derive the desired conclusion from this and (\ref{35.3}), given $\eta>0$ we use (\ref{35.5})
  to find some large $t_0>1$ such that
  \be{35.6}
	\int_{t_0-1}^\infty \|n(\cdot,t)-\mu(t)\|_{L^2(\Omega)}^2 dt <\frac{\eta}{12C_1},
  \ee
  and such that in accordance with Lemma \ref{lem32} we moreover have
  \be{35.7}
	c(x,t) < \frac{\eta}{3C_1}
	\qquad \mbox{for all $x\in\Omega$ and } t>t_0-1
  \ee
  and well as
  \be{35.8}
	\io |\nabla c(\cdot,t)|^2
	<\frac{\eta}{3C_1} 
	\qquad \mbox{for all } t>t_0-1.
  \ee
  Then for arbitrary $t>t_0-1$ we may use (\ref{35.6}) to pick $\ts=\ts(t)\in (t-1,t)\setminus N$ such that
  \be{35.9}
	\io \Big|n(\cdot,\ts)-\mu(\ts)\Big|^2 < \frac{\eta}{12C_1}.
  \ee
  Again by (\ref{24.9}) and the Cauchy-Schwarz inequality, this firstly entails that
  \bas
	|\onz-\mu(\ts)| \cdot |\Omega|
	&=& \bigg| \io \Big(n(\cdot,\ts)-\mu(\ts)\Big) \bigg| \\
	&\le& \bigg\{ \io \Big( n(\cdot,t)-\mu(t)\Big)^2 \bigg\}^\frac{1}{2} \cdot |\Omega|^\frac{1}{2} \\
	&\le& \sqrt{\frac{\eta |\Omega|}{12C_1}}
  \eas
  and hence 
  \bas
	\io \Big(\onz-\mu(\ts)\Big)^2 \le \frac{\eta}{12C_1},
  \eas
  so that, secondly, from (\ref{35.9}) we obtain that
  \bas
	\io \Big(n(\cdot,\ts)-\onz\Big)^2
	&\le& 2 \io \Big(n(\cdot,\ts)-\mu(\ts)\Big)^2
	+ 2\io \Big(\onz -\mu(\ts)\Big)^2 \\[2mm]
	&<& \frac{\eta}{3C_1}.
  \eas
  In conjunction with (\ref{35.7}), (\ref{35.8}) and (\ref{35.3}), this means that
  \bas
	\io \Big(n(\cdot,t)-\onz\Big)^2
	< C_1 \cdot \Big\{ \frac{\eta}{3C_1} + \frac{\eta}{3C_1} + \frac{\eta}{3C_1}\Big\} =\eta,
  \eas
  because $t\in (\ts,\ts+1)\setminus N$.
  Since $\eta>0$ was arbitrary, this completes the proof.
\qed
By interpolation and approximation, in view of the generalized continuity property of $n$ gained in Lemma \ref{lem24}
this readily implies convergence in the style claimed in Theorem \ref{theo99}.
\begin{cor}\label{cor355}
  For all $p\ge 1$,
  \be{355.1}
	n(\cdot,t) \to \onz
	\quad \mbox{in $L^p(\Omega)$ as } t\to\infty.
  \ee
\end{cor}
\proof
  By boundedness of $\Omega$, we only need to consider the case $p>2$, in which due to the H\"older inequality,
  \bas
	\|n(\cdot,t)-\onz\|_{L^p(\Omega)} 
	\le \|n(\cdot,t)-\onz\|_{L^\infty(\Omega)}^\frac{p-2}{p}
	\|n(\cdot,t)-\onz\|_{L^2(\Omega)}^\frac{2}{p}
	\le C_1	
	\|n(\cdot,t)-\onz\|_{L^2(\Omega)}^\frac{2}{p}
	\qquad \mbox{for all } t>0
  \eas
  with $C_1:=\big\{ \|n\|_{L^\infty(\Omega\times (0,\infty))} + \onz \big\}^\frac{p-2}{p}$.
  Therefore, given $\eta>0$ we may invoke Lemma \ref{lem35} to fix $t_0>0$ such that
  \be{355.2}
	\|n(\cdot,t)-\onz\|_{L^p(\Omega)} \le \eta
	\qquad \mbox{for all } t\in (t_0,\infty)\setminus N,
  \ee
  and for the proof of (\ref{355.1}) it will be sufficient to make sure 
  that the inequality herein actually remains valid for all $t>t_0$. 
  To verify this, for any such $t$ we can use the density of $(t_0,\infty)\setminus N$ in $(t_0,\infty)$ to find
  $(t_k)_{k\in\N}\subset (t_0,\infty)\setminus N$ such that $t_k\to t$ as $k\to\infty$. Then (\ref{355.2})
  shows that $\|n(\cdot,t_k)-\onz\|_{L^p(\Omega)} \le \eta$ for all $k\in\N$, whence we may extract a subsequence
  $(t_{k_l})_{l\in\N}$ of $(t_k)_{k\in\N}$ such that $n(\cdot,t_{k_l})-\onz \wto z$ in $L^p(\Omega)$ as $l\to\infty$.
  But since this trivially entails that also $n(\cdot,t_{k_l})-\onz \wto z$ in $(W_0^{2,2}(\Omega))^\star$, from the
  continuity property implied by (\ref{24.3}) we infer that $z$ must coincide with $n(\cdot,t)-\onz$ and that thus
  \bas
	\|n(\cdot,t)-\onz\|_{L^p(\Omega)} \le \liminf_{l\to\infty} \|n(\cdot,t_{k_l})-\onz\|_{L^p(\Omega)} \le \eta,
  \eas
  as claimed.
\qed
\subsection{Decay of $u$}
Finally, uniform decay of $u$ can be achieved on the basis of the following straightforward application of
standard regularity theory in the forced Stokes evolution system.
\begin{lem}\label{lem36}
  There exist $\lambda>0$ and $C>0$ such that for any choice of $\mu\in\R$ and arbitrary $\eps\in (0,1)$ we have
  \be{36.1}
	\|\ueps(\cdot,t)\|_{L^\infty(\Omega)} 
	\le C + C \int_0^t (t-s)^{-\alpha} e^{-\lambda(t-s)} \big\|\neps(\cdot,s)-\mu\big\|_{L^2(\Omega)} ds
	\qquad \mbox{for all } t>0,
  \ee
  where $\alpha \in (\frac{3}{4},1)$ is taken from (\ref{init}).
\end{lem}
\proof
  As gradients of functions from $W^{1,\infty}(\Omega)$ belong to the kernel of the Helmholtz projection,
  for arbitrary $\mu\in\R$ the third equation in (\ref{0eps}) can be rewritten according to
  \be{36.2}
	u_{\eps t} + A\ueps =\proj\Big[ (\neps(\cdot,t)-\mu)\nabla\phi \Big],
	\qquad x\in\Omega, \ t>0.
  \ee
  Now since $\alpha>\frac{3}{4}$,
  from a known embedding result (\cite{giga1981_theother}, \cite{henry}) 
  we obtain that $D(A^\alpha) \hra L^\infty(\Omega)$,
  so that invoking well-known smoothing properties of the analytic semigroup $(e^{-tA})_{t\ge 0}$ (\cite{sohr_book}, 
  \cite{friedman})
  we infer from (\ref{36.2}) that with some $C_1>0$, $C_2>0$ and $\lambda>0$ we have
  \bas
	\|\ueps(\cdot,t)\|_{L^\infty(\Omega)}
	&\le& C_1\|A^\alpha \ueps(\cdot,t)\|_{L^2(\Omega)} \\
	&=& C_1 \bigg\|A^\alpha e^{-tA} u_0
	+ \int_0^t A^\alpha e^{-(t-s)A} \proj\Big[ (\neps(\cdot,s)-\mu)\nabla\phi\Big] ds \bigg\|_{L^2(\Omega)} \\
	&\le& C_1 \|A^\alpha u_0\|_{L^2(\Omega)}
	+ C_2 \int_0^t (t-s)^{-\alpha} e^{-\lambda(t-s)} \Big\| \proj\Big[ (\neps(\cdot,s)-\mu)\nabla\phi\Big] 
		\Big\|_{L^2(\Omega)} ds
  \eas
  for all $t>0$.
  Since $\proj$ is an orthogonal projector and hence 
  \bas
	\Big\| \proj\Big[ (\neps(\cdot,s)-\mu)\nabla\phi\Big] \Big\|_{L^2(\Omega)}
	\le \Big\| (\neps(\cdot,s)-\mu)\nabla\phi \Big\|_{L^2(\Omega)}
	\le \|\nabla\phi\|_{L^\infty(\Omega)} \|\neps(\cdot,s)-\mu\|_{L^2(\Omega)}
	\quad \mbox{for all } s>0,
  \eas
  in view of our regularity assumption $u_0\in D(A^\alpha)$ we thereby obtain (\ref{36.1}).
\qed
Here the integral on the right-hand side can be estimated by using the following elementary decay property.
\begin{lem}\label{lem37}
  Let $\beta\in (0,1)$, $\lambda>0$ and $h:(0,\infty)\to \R$ be measurable and bounded with $h(t)\to 0$ as $t\to\infty$.
  Then
  \be{37.1}
	\int_0^t (t-s)^{-\beta} e^{-\lambda(t-s)} h(s) ds \to 0
	\qquad \mbox{as } t\to\infty.
  \ee
\end{lem}
\proof
  Given $\eta>0$, we pick $t_1>0$ large such that $|h(t)| \le \frac{\eta}{2C_1}$ for all $t>t_1$, where
  $C_1:=\int_0^\infty \sigma^{-\beta} e^{-\lambda\sigma} d\sigma$ is finite since $\beta<1$.
  Then writing $t_0:=t_1+\Big(\frac{2\|h\|_{L^\infty((0,\infty))}}{\lambda}\Big)^\frac{1}{\beta}$,
  for arbitrary $t>t_0$ we can estimate
  \bas
	\bigg| \int_0^t (t-s)^{-\beta} e^{-\lambda(t-s)} h(s) ds \bigg|
	&\le& \int_0^{t_1} (t-s)^{-\beta} e^{-\lambda(t-s)} |h(s)| ds
	+ \int_{t_1}^t (t-s)^{-\beta} e^{-\lambda(t-s)} |h(s)| ds \\
	&\le& (t-t_1)^{-\beta} \|h\|_{L^\infty((0,\infty))} \int_0^{t_1} e^{-\lambda(t-s)} ds
	+ \frac{\eta}{2C_1} \int_{t_1}^t (t-s)^{-\beta} e^{-\lambda(t-s)} ds \\
	&=& (t-t_1)^{-\beta} \|h\|_{L^\infty((0,\infty))} \cdot \frac{1}{\lambda}(1-e^{-\lambda(t-t_1)})
	+ \frac{\eta}{2C_1} \int_0^{t-t_1} \sigma^{-\beta} e^{-\lambda\sigma} d\sigma \\
	&\le& (t_0-t_1)^{-\beta} \|h\|_{L^\infty((0,\infty))} \cdot \frac{1}{\lambda}
	+ \frac{\eta}{2C_1} \int_0^\infty \sigma^{-\beta} e^{-\lambda\sigma} d\sigma \\
	&=& \frac{\eta}{2} + \frac{\eta}{2} =\eta
  \eas
  and thereby see that indeed (\ref{37.1}) is valid.
\qed
In view of the stabilization property from Corollary \ref{cor355}, Lemma \ref{lem36} thus entails the desired
decay feature of $u$.
\begin{lem}\label{lem38}
  We have
  \be{38.1}
	u(\cdot,t)\to 0
	\quad \mbox{in } L^\infty(\Omega)
	\qquad \mbox{as } t\to\infty.
  \ee
\end{lem}
\proof
  With the null set $N\subset (0,\infty)$ taken from Lemma \ref{lem24}, on combining Lemma \ref{lem22}
  with the dominated convergence theorem we obtain that
  $\neps(\cdot,t)-\onz\to n(\cdot,t)-\onz$ in $L^2(\Omega)$ as $\eps=\eps_j\searrow 0$.
  Therefore, using the convergence property (\ref{24.7}) we infer from Lemma \ref{lem36} that there exist
  $\lambda>0$ and $C_1>0$ fulfilling
  \bas
	\|u(\cdot,t)\|_{L^\infty(\Omega)}
	\le C_1 + C_1 \int_0^t (t-s)^{-\alpha} e^{-\lambda(t-s)} \|n(\cdot,s)-\onz\|_{L^2(\Omega)} ds
	\qquad \mbox{for all } t>0,
  \eas
  where $\alpha\in (\frac{3}{4},1)$ is as in (\ref{init}).
  Since $\|n(\cdot,t)-\onz\|_{L^2(\Omega)} \to 0$ as $t\to\infty$ by Corollary \ref{cor355}, Lemma
  \ref{lem37} therefore yields (\ref{38.1}).
\qed
\subsection{Proof of Theorem \ref{theo99}}
We finally only need to collect our previous findings to arrive at our main result.\abs
\proofc of Theorem \ref{theo99}. \quad
  The statement on global existence of a weak solution with the regularity features in (\ref{99.1})
  has been asserted by Lemma \ref{lem24}. The convergence properties in (\ref{99.2}) are precisely 
  established by Corollary \ref{cor355}, Lemma \ref{lem32} and Lemma \ref{lem38}.
\qed
\vspace*{5mm}
{\bf Acknowledgement.} \quad
  The author acknowledges support of the {\em Deutsche Forschungsgemeinschaft} in the context of the project
  {\em Analysis of chemotactic cross-diffusion in complex frameworks}.
\end{document}